\documentclass[11pt]{amsart}

\usepackage{amsmath}
\usepackage{amssymb}
\usepackage{graphicx}
\usepackage{hyperref}
\usepackage{subcaption}
\usepackage[normalem]{ulem}
\usepackage[framemethod=TikZ]{mdframed}
\usepackage{xcolor}
\usepackage[ruled,vlined]{algorithm2e}
\usepackage{enumitem}

\newcommand{\Fix}{\ensuremath{\operatorname{Fix}}}

\usepackage{color}
\definecolor{orange}{rgb}{.8,0.4,0}
\definecolor{purple}{rgb}{.5,0,.5}
\definecolor{violet}{rgb}{.5,0,.5}
\definecolor{OliveGreen}{rgb}{0,0.5,0}
\definecolor{dOliveGreen}{rgb}{0,0.25,0}
\definecolor{dred}{rgb}{0.5,0,0}
\definecolor{dorange}{rgb}{0.5,0.25,0}
\definecolor{dblue}{rgb}{0,0,0.35}
\definecolor{dpurple}{rgb}{0.25,0,0.5}
\definecolor{dmagenta}{rgb}{0.5,0,0.25}

\newtheorem{theorem}{Theorem}[section]

\theoremstyle{definition}
\newtheorem{definition}[theorem]{Definition}
\newtheorem{example}[theorem]{Example}
\newtheorem{remark}[theorem]{Remark}

\numberwithin{equation}{section}

\newcounter{emphasisbox}
\renewcommand{\theemphasisbox}{\thesection.\arabic{emphasisbox}}

\mdfdefinestyle{emphasisbox}{%
	linecolor=blue,
	outerlinewidth=2pt,
	bottomline=false,
	leftline=false,rightline=false,
	skipabove=\baselineskip,
	skipbelow=\baselineskip,
	frametitle=\mbox{},
}
\newmdenv[%
style=emphasisbox,
settings={\global\refstepcounter{emphasisbox}},
frametitlefont={\bfseries emphasisbox~\theemphasisbox\quad},
]{emphasisbox}
\newmdenv[%
style=emphasisbox,
frametitlefont={\bfseries emphasisbox~\quad},
]{emphasisbox*}

\newmdenv[%
backgroundcolor=red!8,
linecolor=red,
outerlinewidth=1pt,
roundcorner=5mm,
skipabove=\baselineskip,
skipbelow=\baselineskip,
]{emphasisboxed}

\title[Projections and differential equations]{Application of projection  algorithms to differential equations:\\boundary value problems}

\author[Bishnu P. Lamichhane]{Bishnu P. Lamichhane}
\address[Bishnu P. Lamichhane]{CARMA, University of Newcastle, Australia}
\email{{\tt bishnu.lamichhane@newcastle.edu.au}}

\author[Scott B. Lindstrom]{Scott B. Lindstrom}
\address[Scott B. Lindstrom]{CARMA, University of Newcastle, Australia}
\email{{\tt scott.lindstrom@uon.edu.au}}

\author[Brailey Sims]{Brailey Sims}
\address[Brailey Sims]{CARMA, University of Newcastle, Australia}
\email{\tt brailey.sims@newcastle.edu.au}



\date{\today}

\def\R{\hbox{$\mathbb R$}}

\newcommand{\cL}{\mathcal{L}}

\newcommand{\cS}{\mathcal S}

\begin{document}
	
	\maketitle

	\begin{abstract}
	The Douglas-Rachford method has been employed successfully to solve many kinds of non-convex feasibility problems. In particular, recent research has shown surprising stability for the method when it is applied to finding the intersections of hypersurfaces. Motivated by these discoveries, we reformulate a second order boundary value problem (BVP) as a feasibility problem where the sets are hypersurfaces. We show that such a problem may always be reformulated as a feasibility problem on no more than three sets and is well-suited to parallelization. We explore the stability of the method by applying it to several examples of BVPs, including cases where the traditional Newton's method fails.
	\end{abstract}

	\section{Introduction}

We explore a particular approach to obtaining approximate numerical solutions to (second order, nonlinear) boundary-value problems on $[a,b]\subseteq \mathbb{R}$. We use finite difference approximations to replace the continuous problem by a discrete one involving a  finite system of $N$ nonlinear equations in $N$ variables (the approximate solution values at each of the $N$ partition points). The classical approach to solving such a system of equations is to use Newton's method. We explore some alternative projection-based iterative methods. 

The solution set for each of the $N$ equations is a hypersurface $S_k$ in $N$-dimensional Euclidean space $\mathbb{R}^N$. An approximate solution to the BVP then corresponds to a point in the intersection of these $N$ hypersurfaces. This is a \emph{feasibility} problem of the form: 
\begin{equation*}
	\text{Find}\quad x \in \bigcap_{k=1}^N S_k.
\end{equation*}
One approach to solving such feasibility problems is to use an iterated process. 

We consider the method of alternating projections (AP) and the Douglas-Rachford method (DR) in particular. We explain how the methods are well-suited to parallelization. We then use the methods to solve the associated feasibility problems for several BVPs and compare the results with those given by the classical Newton's method. 

\subsection{Objectives}

Our intent is not to compare the speeds of our projection-based methods with that of Newton's method, which is much faster. Neither is it our goal to provide a full comparison of their respective robustness. The main contributions are as follows.
\begin{enumerate}
	\item We introduce the reformulation of ODEs as hypersurface feasibility problems for solving with iterated projection methods.
	\item We show how they are particularly amenable to parallelization.
	\item We show how, for boundary value ODE's, we may reformulate the $N$ hypersurface feasibility problem as a $3$-set feasibility problem
	\item We analyse the behaviour for both AP and DR experimentally on hypersurface problems for varying $N$, which for boundary value ODEs corresponds to partition fineness. We catalogue the characteristics of oscillation so frequently observed for DR in particular.
	\item We provide a characterization of how it might be employed to real world systems of equations in cases where Newton's method does not succeed.
\end{enumerate}
This work extends to $N$ sets --- via Pierra's method \cite{Pierra} (aka the \emph{divide and concur} method \cite{GE}) --- the $2$ set investigation started by Borwein and Sims \cite{BS}, who analysed DR for the hypersurfaces choices of a $(n-1)$-sphere and a line. In this simpler setting, global convergence was shown by Borwein and Aragon Artacho \cite{AB} under an assumption later relaxed by Benoist \cite{Benoist}, who demonstrated convergence by means of a Lyapunov function. The analysis has already been extended in $\mathbb{R}^2$ by Borwein, Lindstrom, Schneider, Sims, and Skerritt \cite{BLSSS}, who considered the generalization of circles to ellipses and $p$-spheres. Later, Lindstrom, Sims, and Skerritt considered plane curves more generally \cite{LSS}. Inspired by Benoist's work, Dao and Tam \cite{DT} have since provided a beautiful illumination of the method for curves in $\mathbb{R}^2$ by means of Lyapunov functions. Phan \cite{Phan} and later Phan and Dao \cite{DP} have since provided more general convergence results under regularity (transversality) conditions. 

While this article is an important extension of the analysis of projection methods (and DR in particular) for nonconvex hypersurface feasibility problems, its approach is comparable to other experimental works which analyse proximal point algorithms in the absence of convexity by cataloguing the performance of the method for a selection of examples. These include the recent work of Arag\'on Artacho, Borwein, and Tam applying DR to solve matrix completion problems \cite{ABT2} and Sudoku puzzles \cite{ABT1}, the work of Arag\'on Artacho and Campoy with graph coloring problems \cite{AC}, and the seminal work of Elser, Rankenburg, and Thibault \cite{ERT}. 

We have listed here only a small selection of the nonconvex Douglas-Rachford genre. The history is vast, and we have not even touched on its roots in convex optimization and connections with the celebrated ADMM through duality. For a more thorough treatment, we refer the reader to a recent survey of Sims and Lindstrom \cite{LS}.

\subsection{Outline}

The outline of this paper is as follows. In \ref{s:nonlinearboundaryvalueproblems} we introduce nonlinear boundary value problems. In \ref{s:projectionmethods} we introduce the $2$ set projection algorithms and their extension to $N$ sets. In \ref{s:computationsfornearestpoint}, we discuss methods of projecting onto individual hypersurfaces, and in \ref{s:theprocedure} we describe the full procedure, discuss its amenability to parallelization, and show a natural reformulation which reduces the $N$ set problem to a $3$ set problem. We provide our experimental results and conclude in \ref{s:examples}.


\section{Nonlinear boundary value problems}\label{s:nonlinearboundaryvalueproblems}

We investigate the use of projection algorithms to obtain numerical solutions to
nonlinear boundary value problems. Here and throughout:
\begin{emphasisboxed}
	\begin{enumerate}[label={(\roman*)}]
		\item $y:[a,b]\subset \mathbb{R}\rightarrow\mathbb{R}$ with $a<b$ is an ``unknown'' function for which we seek a numerical solution.
		\item $y'$ and $y''$ are, respectively, the first and second derivatives of $y$.
		\item $\alpha:= y(a) \in \mathbb{R}$ and $\beta:=y(b) \in \mathbb{R}$ are given boundary values.
	\end{enumerate}
\end{emphasisboxed}
A complete statement of the problem is:
\begin{emphasisboxed}
	\begin{align}
		\text{Find}\quad &y \;{\rm such \; that} \; \nonumber\\
		&y'' = f(x,y,y')\quad\text{for}\;  x \in (a,b) \subset \mathbb{R} \; \text{with}\; y(a) = \alpha\;\text{and}\; y(b) =\beta.\label{bvp}
	\end{align}
\end{emphasisboxed}

\begin{remark}[Solutions may not be unique]\label{rem:uniqueness}
	In general, even when a solution to problem \eqref{bvp} exists, it may not be unique. However, \eqref{bvp} will have a unique continuous solution over the interval $[a,b]$ if the right-hand side function $f$ satisfies the following conditions:
	\begin{enumerate}
		\item $f$ and the partial derivatives of $f$ with respect to $y$ and $y'$ are continuous on
		\[ D = \{(x,y,y')\,|\, a\leq x \leq b, \; -\infty < y <\infty,\;  -\infty < y' <\infty\},\]
		\item $\frac{\partial f}{\partial y}(x,y,y') >0$ on $D$, and
		\item there exists a constant $M$ such that
		\[ \left|\frac{\partial f}{\partial y'}(x,y,y') \right| \leq M \quad\text{on}\; D. \]
	\end{enumerate}
	See, for example, \cite[Theorem 11.1]{BF16}. Because we seek to present the wide variety of behaviours exhibited by our algorithms, we will present both examples which do and do not satisfy these criteria.
\end{remark}
We use a finite difference method to approximate the solution of the given boundary value problem.
This results in  a system of nonlinear equations to which we apply our projection algorithm to compute an approximate numerical solution.

To this end, consider a partition of the interval $[a,b]$ into $N$ equal subintervals using the set of points
$x_i = a+ih$ for $i=0,1,\cdots,N+1$ with $x_{N+1} =b$ so that
\[ h = \frac{b-a}{N}.\]

We introduce the centred-difference approximations,
$$
y'(x_i) \approx \frac{y(x_{i+1}) - y(x_{i-1})}{2h}
$$
and
$$
y''(x_i) \approx \frac{y(x_{i+1}) - 2y(x_i) + y(x_{i-1})}{h^2}.
$$
When the exact solution $y$ is four times continuously differentiable these estimate the first and second derivatives at $x_{i}$ with errors of $\frac{h^2}{6} y^{(3)}(\eta_i)$ and $\frac{h^2}{12} y^{(4)}(\xi_i)$ respectively, where $\eta_i$ and $\xi_i$ lie in the interval $(x_{i-1},x_{i+1})$.

Ignoring such truncation error terms we replace the first and second derivative of $y$ by their centred-difference approximations in \eqref{bvp} to obtain for $i =1,2,3,\cdots, N$ the approximate relationships
\[ \frac{y(x_{i+1}) - 2y(x_i) + y(x_{i-1})}{h^2} \approx f\left(x_i, y(x_i) ,  \frac{y(x_{i+1}) - y(x_{i-1})}{2h} \right).\]
This leads us to take as an approximate numerical solution to \eqref{bvp} $y(x_{i}) \approx \omega_{i}$ where the $\omega_{i}$ satisfy the system of generally nonlinear equations

$\omega_0 = \alpha$, $\omega_{N+1} = \beta$ and
\begin{equation}
	\frac{\omega_{i+1} -2 \omega_i + \omega_{i-1}}{h^2} - f\left(x_i, \omega_i, \frac{\omega_{i+1} - \omega_{i-1}}{2h}\right) =0  \quad\text{for}\;
	i =1,2,3,\cdots, N.
\end{equation}\label{nbvp}
If $h<2/M$ where $M$ is as defined in Remark~\ref{rem:uniqueness} and the other conditions of Remark~\ref{rem:uniqueness} are satisfied, then this nonlinear system of equations has a unique solution \cite[page 86]{Keller68}. While many of our examples do not satisfy the conditions of Remark~\ref{rem:uniqueness}, uniqueness implies that we can easily measure the accuracy of a numerical approach, whereas if we have non-uniqueness it is much harder.
\begin{align}
	\varphi_i(\omega)&:=\frac{\omega_{i+1} -2 \omega_i + \omega_{i-1}}{h^2} - f\left(x_i, \omega_i, \frac{\omega_{i+1} - \omega_{i-1}}{2h}\right),\label{sys} \\
	&\textrm{for\;} i = 1,2,\cdots, N,\nonumber\\
	&\textrm{where}\;\omega_0=\alpha,\omega_{N+1}=\beta,\; \textrm{and set\;}\nonumber \\
	\Omega_i &:= \{\omega=(\omega_1,\dots,\omega_N)\;|\;\omega \;{\rm satisfies\;}\; \varphi_i(\omega) = 0\}, \label{feasibility_reformulation}\\
	&\textrm{for\;} i = 1,2,\cdots, N. \nonumber
\end{align}
Then we can compute our approximate numerical solution to the boundary value problem  \eqref{bvp} by solving the feasibility problem: find $\omega \in \cap_{i=1}^N \Omega_i$. An approximate numerical solution to \eqref{bvp} is then given by $y(x_i)={\omega}_{i}$. For the task, we employ both the method of alternating projections and a parallelized version of the Douglas-Rachford method as outlined below.

\begin{remark} The astute reader will note that more complicated boundary conditions may be handled by appropriately modifying either or both of the equations $\omega_0 = \alpha$, $\omega_{N+1}=\beta$ though this could potentially lead to an enlarged problem of $N+2$ equations in $N+2$ unknowns. For example, the mixed condition $y(a) + \eta y'(a) = \alpha$ could translate to $\displaystyle \varphi_0(\omega) = \omega_{0} + \frac{\omega_{1} - \omega_{0}}{h}   = \alpha$.
\end{remark}

\section{Preliminaries on Projection Methods}\label{s:projectionmethods}

The Douglas-Rachford method (DR) and the method of alternating projections (AP) are frequently used to find a feasible point (point in the intersection) of two closed constraint sets $A$ and $B$ in a Hilbert space, in our setting: $N$-dimensional Euclidean space, $\mathbb{R}^N$.

The projection onto a subset $C$ of $\mathbb{R}^N$ is defined for all $x \in \mathbb{R}^N$ by $$\mathbb{P}_C(x) := \left \{ z \in C : \|x - z\| = \inf_{z' \in C}\|x - z'\|\right \}.$$

\begin{emphasisboxed}
	Note that $\mathbb{P}_C$,is a set-valued map where values may be empty or contain more than one point. In our case of interest, where $C$ is a closed hypersurface, $\mathbb{P}_C$ has nonempty values and, in order to simplify both notation and implementation, we will work with a selector for $\mathbb{P}_C$ that is a map $P_C:\mathbb{R}^N \rightarrow C: x \mapsto P_C(x) \in \mathbb{P}_C(x)$.
\end{emphasisboxed}

When $C$ is nonempty, closed, and convex the projection operator $P_C$ is uniquely determined and \emph{firmly nonexpansive}; that is $\left(\forall x,y\in \mathbb{R}^N \right)$
\begin{equation*}
	\|P_C x - P_C y \|^2 + \|(I-P_C)x-(I-P_C)y\|^2 \leq \|x-y\|^2.
\end{equation*}
See, for example, \cite[Chapter 4]{BC}. When $C$ is a closed subspace it is also a linear operator \cite[Corollary 3.22]{BC}.

The reflection mapping through the set $C$ is then defined by
$$R_C := 2P_C - I,$$
where $I$ is the identity map.

\begin{definition}[Method of Alternating Projections]
	For two closed sets $A$ and $B$ and an initial point $x_0 \in H$, the method of alternating projections (AP) generates a sequence $(x_n)_{n=1}^\infty$ as follows:
	\begin{equation}
		x_{n+1} \in T'_{A,B}(x_n) \quad \text{where} \quad T'_{A,B} := P_{B}P_{A}.
	\end{equation}	
\end{definition}

The Douglas-Rachford method was introduced half a century ago in connection with nonlinear heat flow problems \cite{DR}.

\begin{definition}[Douglas-Rachford Method]
	For two closed sets $A$ and $B$ and an initial point $x_0 \in H$, the Douglas-Rachford method (DR) generates a sequence $(x_n)_{n=1}^\infty$ as follows:
	\begin{equation}
		x_{n+1} \in T_{A,B}(x_n) \quad \text{where} \quad T_{A,B} := \frac{1}{2}\left( I + R_{B}R_{A}\right).
	\end{equation}	
\end{definition}

\begin{definition}[Fixed point set]
	The fixed point set for an operator $T$ is $\Fix T = \{x \in H | x \in Tx \}$.
\end{definition}

The following theorem from \cite{BCL} relaxes, in the context of convex feasibility, previous convergence conditions established in the somewhat different context of \cite{LM}. See also \cite{svaiter}.

\begin{theorem}[{\cite[Fact 5.9]{BCL}}]
	Suppose $A,B \subseteq H$ are closed and convex with non-empty intersection. Given $x_0 \in H$ the sequence of iterates defined by $x_{n+1}:=T_{A,B}x_n$ converges weakly to an $x \in \Fix T_{A,B}$ with $P_{A}x \in A \cap B$.
\end{theorem}
Of course in our case, where the space is finite dimensional, weak convergence ensures convergence in norm.

Notwithstanding the absence of a satisfactory theoretical justification, the Douglas-Rachford iteration scheme has been used to successfully solve a wide variety of practical problems in which one or both of the constraints are non-convex.

In an effort to develop the beginnings of a theoretic basis for employment in the non-convex setting, the authors of \cite{BS} explored a feasibility problem on two particular hypersurfaces in $\mathbb{R}^n$:  a line and the $n-1$-sphere. Among other results, they established local convergence near each of the (possibly two) feasible points. More extensive regions of convergence were determined by Borwein and Arag\'{o}n Artacho \cite{AB}. The definitive answer, as conjectured in \cite{BS}, was subsequently  given by Benoist \cite{Benoist} who established convergence to the nearest feasible point except for starting points lying on a singular set: the hyperplane of symmetry. 

Borwein et al.\ \cite{BLSSS} showed that local convergence still holds for a line and a smooth  hypersurface in $\mathbb{R}^N$ not intersecting asymptotically, although the basins of convergence may be quite sensitive to small perturbations of the sets. Additionally, Lindstrom et al.\ \cite{LSS} extended local convergence to isolated points of intersection for two smooth hypersurfaces in $\R^{N}$. The authors of \cite{LLM} showed local convergence for the John von Neumann’s method of alternating projections for sets under regularity conditions. Phan \cite{Phan}, and later Phan and Dao \cite{DP}, showed local convergence with \emph{R-linear} convergence rate for the strongly regular system $\{A,B\}$ of superregular sets $A,B$. For more details on the history, we again refer the reader to \cite{LS}.

\subsection{Extension of DR to Many Sets}
We can apply this method to a consistent feasibility problem with $N$ sets $\Omega_1 \dots \Omega_N\subset \mathbb{R}^N$ to find $x \in \cap_{k=1}^N \Omega_k \ne \emptyset$. We do so by working in the product space $\mathbb{R}^{N\times N}$ as follows.
\begin{align}
	\text{Let}\quad	A&:=\Omega_1 \times \dots \times \Omega_N\nonumber \\
	\text{and}\quad	B&:= \{(x_1,\dots,x_N)\in \mathbb{R}^{N\times N} | x_1=x_2=\dots = x_N \}\label{eqn:B}
\end{align}
and apply the DR method to the two sets $A$ and $B$. The product space projections for $x=(x_1,\dots,x_N)\in  \mathbb{R}^{N\times N}$ are
\begin{align*}  P_A(x_1,\dots,x_N)&=(P_{\Omega_1}(x_1),\dots, P_{\Omega_N}(x_N))\\
	P_B(x_1,\dots,x_N)&=\left(\frac{1}{N}\sum_{k=1}^N x_k,\dots,\frac{1}{N}\sum_{k=1}^N x_k\right).
\end{align*}
This is sometimes called the ``divide and concur'' method. See, for example, \cite{Pierra,GE}. The ``divide and concur'' method is particularly well suited to parallelization. An alternative would be to use the cyclic Douglas-Rachford algorithm introduced in \cite{BTam}.

We consider, in particular, the case where the $\Omega_i$ are as in \eqref{feasibility_reformulation}. Where $\omega_0=\alpha$ and $\omega_{N+1}=\beta$ are fixed, the \textsl{feasibility problem} is reduced to finding a point in the intersection of a family of $N$ hyper-surfaces $\Omega_1,\dots,\Omega_N$ in $\R^{N}$.

\section{Computations for nearest point projection onto a surface $P_{\Omega_k}$}\label{s:computationsfornearestpoint}

For a hypersurface $\Omega_k$ in $X=\R^{N}$ implicitly defined by 
\begin{equation*}
	x \in \Omega_k \iff \phi_k(x)=0
\end{equation*}
the nearest point projection $x=P_{\Omega_k}(u)$ of $u\in \R^{N}$ onto $\Omega_k$ solves
\begin{align}
	\textrm{minimize:}&\quad F(x)=\frac{1}{2}\|u - x\|^{2} \nonumber\\
	\textrm{subject to:}&\quad \phi_k(x)\ =\ 0.\label{subproblem}
\end{align}
So, provided $u \not\in \Omega_k$ and assuming sufficient differentiability, we know by the theory of Lagrange multipliers --- see, for example, \cite{RV1973} --- that there exists  $\lambda_{0} \not=0$ for which $(P_{\Omega_k}(u),\, \lambda_{0})$ is a critical point of the Lagrangian $\cL(x, \lambda) := F(x)-\lambda\phi_k(x)$. That is,  $\lambda = \lambda_{0}$ and $x= P_{\Omega_k}(u)$ is a solution of,
\begin{equation}
	u - x + \lambda\nabla_{x}\phi_k(x)\ =\ 0,\ \textrm{and\ }\phi_k(x)\ =\ 0.\label{LMult}
\end{equation}


Again assuming sufficient differentiability to ensure that the $(N+1)\times(N+1)$ Jacobian
\begin{align*}
	J(x,\lambda)_{l,j} = \left( \begin{array}{cc} \lambda\left[\frac{\partial^{2}}{\partial x_{l}\partial x_{j}}\phi \right] - I& \left[\nabla_{x}\phi_l(x)\right]^{T} \\ \nabla\phi_{x}(x)& 0\end{array}\right)
\end{align*}
is well-defined, the nonlinear system \eqref{LMult} could be solved using Newton's method. However, this requires solving --- at each iteration --- the system of $N+1$ equations given by $J(x,\lambda)v=b$. 
Quasi-Newton method requires solving a similar linear system.

Alternatively we might seek to locate a point $(x_{0},\lambda_{0})$ where the scalar function
$$
G(x,\lambda)\ :=\ \frac{1}{2}\left(\left\|u - x + \lambda\nabla_{x}\phi(x)\right\|^{2}\ +\ \phi(x)^{2}\right)
$$
has a minimum zero; $P_{\cS}(u) = x_{0}$ is then the desired solution. For this we might use the method of gradient (steepest) descent with a line search implemented at each iteration. This obviates the need to invert $J(x,\lambda)$, but depending on the method employed for the line search may involve performing several iterations of Newton's method on a one variable function at each step. The main difficulty here is choosing a suitable starting point; $(u, 0)$ is one choice.

A simple code for computing the hypersurface projections for $\Omega_1,\dots,\Omega_{N}$ may be seen in Algorithm~\ref{alg:generatehypersurfaces}.

\begin{algorithm}\caption{Compute $\Omega_1,\dots,\Omega_N \subset \mathbb{R}^N$ and projections onto them}\label{alg:generatehypersurfaces}
	\SetKwFor{proc}{procedure}{}{end:}
	\SetKwFor{uFor}{for}{do}{end:}
	\proc{Generate Hypersurfaces}{
		\KwData{Receives as input a function $f$ which defines the differential equation, boundary points $a,b$ with corresponding solution values $\alpha,\beta$ which define the boundary conditions, a number $N$ of partition points.}
		\KwResult{Returns $\phi=(\phi_1,\dots,\phi_{N})$ where $\Omega_k = \{z | \phi_k(z)=0\}$ is the $k$th hypersurface for the feasibility problem.}
		set $h:=\frac{b-a}{N+1}$;\\
		\uFor{$k \in \{1,\dots,N\}$ }{
			\uIf{$k=1$}{
				set $\phi_k:=x\mapsto 2x_k -x_{k+1}+h^2f\left( a+kh,x_k,\frac{x_{k+1}-\alpha}{2h}\right) -\alpha$;
			}
			\uElseIf{$k=N$}{
				set $\phi_k:=x\mapsto 2x_k -\beta+h^2f\left( a+kh,x_k,\frac{\beta-x_{k-1}}{2h}\right) -x_{k-1}$;
			}
			\uElse{
				set $\phi_k:=x\mapsto 2x_k - x_{k+1} + h^2 f\left(a+kh, x_k, \frac{x_{k+1}-x_{k-1}}{2h} \right)- x_{k-1}$;
			}	
		}	
		store $\phi$;
	}
	\proc{Compute Lagrangian Problems}{
		\KwData{receives as input the $N$-tuple $\phi$}
		\KwResult{Stores $\varphi=(\varphi_1,\dots,\varphi_{N})$ where $\varphi_k = \{\varphi_{k,1},\varphi_{k,2},\varphi_{k,3}\}$ is three of the four functions from the Lagrangian system for computing a projection onto $\Omega_k$ (the fourth function is $\phi_k$).}
		\uFor{$k \in \{1,\dots,N\}$ }{
			set $\varphi_{k,2}:=(v,u,\lambda) \mapsto 2u_{k}-2v_{k}-\lambda \partial_{k}\phi_k (u)$;\\
			\uIf{$k=1$}{
				set $\varphi_{k,1}:=(v,u,\lambda) \mapsto 0$;
			}
			\uElseIf{$k=N$}{
				$\varphi_{k,3}:=(v,u,\lambda) \mapsto 0$;
			}
			\uElse{
				set $\varphi_{k,1}:=(v,u,\lambda) \mapsto 2u_{k-1}-2v_{k-1}-\lambda \partial_{k-1}\phi_k (u)$;\\
				set $\varphi_{k,3}:=(v,u,\lambda) \mapsto 2u_{k+1}-2v_{k+1}-\lambda\partial_{k+1}\phi_k (u)$;	
			}
		}
		store $\varphi$;
	}
	\proc{Projection for $\Omega_k$}{
		\KwData{receives as input a value $k \in \{1,\dots,N\}$ and a value $x \in \mathbb{R}^N$.}
		\KwResult{Returns a point $u \in P_{\Omega_k}(x)$.}
		Numerically solve the system $\{\phi_k(u)=0, \varphi_{k,1}(x,u,\lambda)=0,\varphi_{k,2}(x,u,\lambda)=0,\varphi_{k,3}(x,u,\lambda)=0\} $ for $u$;\\
		One may use, for example, Algorithm~\ref{alg:hypersurfaceproject} or Algorithm~\ref{alg:hypersurfaceproject2}.\\
		return $u$;
	}
\end{algorithm}

\begin{algorithm}\caption{Projects a point $x$ onto a set $\Omega_k$ with Newton's method}\label{alg:hypersurfaceproject}
	\SetKwFor{proc}{procedure}{}{end:}
	\SetKwFor{uFor}{for}{do}{end:}
	\SetKwFor{uWhile}{while}{do}{end:}
	\proc{Projection for $\Omega_k$}{
		\KwData{receives as input a value $k \in \{1,\dots,N\}$, a value $x \in \mathbb{R}^N$, and a threshold $\Gamma$}
		\KwResult{Returns a point $u \in P_{\Omega_k}(x)$.}
		Set $G:=(v,\lambda)\mapsto \left(\phi_k(v), \varphi_{k,1}(x,v,\lambda),\varphi_{k,2}(x,v,\lambda),\varphi_{k,3}(x,v,\lambda)\right)$.\\
		Set $J:= (v,\lambda) \rightarrow J(v,\lambda)$ where $J(v,\lambda)$ is the Jacobian of $G$ evaluated at $(v,\lambda)$.\\
		Set ${\rm Newt}:=(v,\lambda)\mapsto (v,\lambda)-\Big(\big({\rm MatrixInverse}(J(v,\lambda))\big).G(v,\lambda)\Big)$ where the dot denotes multiplication of a vector by a matrix.\\
		Set $\eta_{old} = (x,1)$.\\
		Set $\eta_{new} = {\rm Newt}(x,1)$.\\
		\uWhile{$\|\eta_{old}-\eta_{new}\| > \Gamma$}{
			Set $\eta_{old}:=\eta_{new}$.\\
			Set $\eta_{new}:= {\rm Newt}(\eta_{old})$.
		}
		
		Set $u:=\eta_{new}$.\\
		Return $(u_1,\dots,u_N)$, the first $N$ components of $u$. Note that the $(N+1)$th component was merely the final Lagrange multiplier)	
	}
\end{algorithm}

\begin{algorithm}\caption{Projects a point $x$ onto a set $\Omega_k$ with steepest descent method}\label{alg:hypersurfaceproject2}
	\SetKwFor{proc}{procedure}{}{end:}
	\SetKwFor{uFor}{for}{do}{end:}
	\SetKwFor{uWhile}{while}{do}{end:}
	\proc{Projection for $\Omega_k$}{
		\KwData{receives as input a value $k \in \{1,\dots,N\}$, a value $x \in \mathbb{R}^N$, a step size modifier $\gamma$, and a threshold $\Gamma$}
		\KwResult{Returns a point $u \in P_{\Omega_k}(x)$.}
		
		Set $G:=(v,\lambda)\mapsto \left(\phi_k(v)\right)^2+\left( \varphi_{k,1}(x,v,\lambda)\right)^2+\left(\varphi_{k,2}(x,v,\lambda)\right)^2+\left(\varphi_{k,3}(x,v,\lambda)\right)^2$.\\
		
		Set $G'(v,\lambda)$ as the gradient of $G$ evaluated at $(v,\lambda)$.
		
		Set ${\rm Descent}:= (v,\lambda) \mapsto v - \gamma G'(v,\lambda)$.\\
		
		Set $\eta_{old} = (x,1)$.\\
		Set $\eta_{new} = {\rm Descent}(x,1)$.\\
		\uWhile{$\|\eta_{old}-\eta_{new}\| > \Gamma$}{
			Set $\eta_{old}:=\eta_{new}$.\\
			Set $\eta_{new}:= {\rm Descent}(\eta_{old})$.
		}
		
		Set $u:=\eta_{new}$.\\
		Return $(u_1,\dots,u_N)$, the first $N$ components of $u$. Note that the $(N+1)$th component was merely the final Lagrange multiplier)

	}

\end{algorithm}

\section{The procedure}\label{s:theprocedure}
To move from a given iterate to a successive iterate, one must compute the approximate projections $P_{\Omega_k}(u), k=1,\dots,N$. One may use an appropriate iterative numerical method to solve the subproblem \eqref{subproblem}, continuing the method until successive iterates differ by less than some pre-prescribed tolerance $\tau$.

The choice of numerical method is between needing more steps but less computational complexity (without Jacobian) versus needing fewer steps with each entailing greater computational complexity (with the Jacobian). For the sake of simplicity, we used the Jacobian for all of our experiments and computed until the change from one step to the next was less than $10^{-12}$.

Even though $\mathcal{P}_{\Omega_{k}}(x_{m,k})$, as a (possibly rough) numerical approximation to $P_{\Omega_{k}}$, may not lie exactly on the surface $\Omega_{k}$, we naturally use $\mathcal{R}(x_{m,k}):=\left(2\mathcal{P}_{\Omega_{k}} - I\right)(x_{m,k})$ in place of the reflection of $x_{m,k}$  in $\Omega_{m,k}$ when computing the $(m+1)^{\textsl{th}}$ iterate of the Douglas--Rachford algorithm, so that
$$
x_{m+1,k} =\ \frac{1}{N}\left(\sum_{j=1}^{N}\mathcal{R}_{\Omega_{j}} (x_{m,k}) \right) - \frac{1}{2}\mathcal{R}_{\Omega_{k}}\left(x_{m,k}\right) + \frac{1}{2}x_{m,k}.
$$

\begin{remark} One might consider using a tolerance $\tau_{m}$ that reduces as the number of iterations increases and hopefully move nearer to a solution. Otherwise, it's is unlikely that the accuracy of the solution found would exceed the preselected tolerance, $\tau$. One could use $\tau_m=\alpha\textrm{diam}\{x_{m,k}:k=1,2,\cdots,N\}$ where $\alpha\in(0,1)$ and
	\begin{equation}
		\textrm{diam}{S}=\underset{s_i,s_j \in S}{\max}\|s_i-s_j\|.
	\end{equation}
\end{remark}
While theory does not guarantee convergence with either method of computing projections, experimentation has shown that for some of the problems Douglas-Rachford method may be relatively insensitive to small changes in how projections are computed \cite{LSS}. This is why it makes sense to consider adapting the tolerance over successive iterates.

\subsection{Alternative Formulation}

We may also attempt to speed up convergence by considering two modified versions of the method. Consider the problem with a partition of $7$ segments, so $N=6$. From the form of equation~\eqref{sys}, for a single iteration, the values updated by an iteration $x\mapsto R_A(x)$ in the product space are underlined in the table below.
\begin{center}
	\begin{tabular}{l l l l l l l}
		\textcolor{red}{$P_{\Omega_1}(x_1)$}&\textcolor{orange}{$P_{\Omega_2}(x_2)$}&\textcolor{OliveGreen}{$P_{\Omega_3}(x_3)$}&\textcolor{blue}{$P_{\Omega_4}(x_4)$}&\textcolor{purple}{$P_{\Omega_5}(x_5)$}&\textcolor{magenta}{$P_{\Omega_6}(x_6)$}\\ \hline
		\textcolor{red}{\underline{$x_{1_1}$}}&\textcolor{orange}{\underline{$x_{2_1}$}}&$x_{3_1}$&$x_{4_1}$&$x_{5_1}$&$x_{6_1}$\\
		\textcolor{red}{\underline{$x_{1_2}$}}&\textcolor{orange}{\underline{$x_{2_2}$}}&\textcolor{OliveGreen}{\underline{$x_{3_2}$}}&$x_{4_2}$&$x_{5_2}$&$x_{6_2}$\\
		$x_{1_3}$&\textcolor{orange}{\underline{$x_{2_3}$}}&\textcolor{OliveGreen}{\underline{$x_{3_3}$}}&\textcolor{blue}{\underline{\underline{$x_{4_3}$}}}&$x_{5_3}$&$x_{6_3}$\\
		$x_{1_4}$&$x_{2_4}$&\textcolor{OliveGreen}{\underline{$x_{3_4}$}}&\textcolor{blue}{\underline{\underline{$x_{4_4}$}}}&\textcolor{purple}{\underline{\underline{$x_{5_4}$}}}&$x_{6_4}$\\
		$x_{1_5}$&$x_{2_5}$&$x_{3_5}$&\textcolor{blue}{\underline{\underline{$x_{4_5}$}}}&\textcolor{purple}{\underline{\underline{$x_{5_5}$}}}&\textcolor{magenta}{\underline{\underline{$x_{6_5}$}}}\\
		$x_{1_6}$&$x_{2_6}$&$x_{3_6}$&$x_{4_6}$&\textcolor{purple}{\underline{\underline{$x_{5_6}$}}}&\textcolor{magenta}{\underline{\underline{$x_{6_6}$}}}
\end{tabular}\end{center}
However, in the computation of the projection onto the agreement space ($P_B$) values are averaged across rows, and so many unchanged values are included in the averaging step. More precise solutions require higher $N$, and for higher $N$ the ratio of unchanged values to changed values in the averaging step grows. This usually slows down computation and convergence. One possible solution is to reformulate the problem as a problem of computation with three sets, $\Omega_1 \cap \Omega_4 , \Omega_2 \cap \Omega_5$, and $\Omega_3 \cap \Omega_6$, as detailed below.
\begin{center}
	\begin{tabular}{c c c}
		$P_{\textcolor{red}{\Omega_1}\cap \textcolor{blue}{\Omega_4}}(x_1)$ & $P_{\textcolor{orange}{\Omega_2}\cap \textcolor{violet}{\Omega_5}}(x_2)$ & $P_{\textcolor{OliveGreen}{\Omega_3}\cap \textcolor{magenta}{\Omega_6}}(x_3)$\\ \hline
		\textcolor{red}{\underline{$x_{1_1}$}} & \textcolor{orange}{\underline{$x_{2_1}$}} & {$x_{3_1}$} \\
		\textcolor{red}{\underline{$x_{1_2}$}} & \textcolor{orange}{\underline{$x_{2_2}$}} & \textcolor{OliveGreen}{\underline{$x_{3_2}$}} \\
		\textcolor{blue}{\underline{\underline{$x_{1_3}$}}} & \textcolor{orange}{\underline{$x_{2_3}$}} & \textcolor{OliveGreen}{\underline{$x_{3_3}$}} \\
		\textcolor{blue}{\underline{\underline{$x_{1_4}$}}} & \textcolor{violet}{\underline{\underline{$x_{2_4}$}}} & \textcolor{OliveGreen}{\underline{$x_{3_4}$}} \\
		\textcolor{blue}{\underline{\underline{$x_{1_5}$}}} & \textcolor{violet}{\underline{\underline{$x_{2_5}$}}} & \textcolor{magenta}{\underline{\underline{$x_{3_5}$}}} \\
		{$x_{1_6}$} & \textcolor{violet}{\underline{\underline{$x_{2_6}$}}} & \textcolor{magenta}{\underline{\underline{$x_{3_6}$}}}\\
\end{tabular}\end{center}
Here the updated values in each column which are underlined twice may be calculated separately from those underlined once, and so this reformulation is no less amenable to parallelization. Notice that we can reformulate in this way for any $N>3$, and that still only two unchanged values will remain at each step (one for the first partition point and one for the last). The memory necessary to store this product space vector $x$ is smaller, although the number of projections computed remains the same because the computation of $P_{\cap_k \Omega_{i+3k} }(x_i)$ requires the computation of  $P_{\Omega_{i}}(x_i), P_{\Omega_{i+3}}(x_i), \dots $.

Another approach is to simply change the map $P_B$ so that only the changed row values are averaged in the agreement step. $P_B$ is, in this reformulation, still a map into $B$. It is no longer the projection map, but we expect similar behavior to that of the three set reformulation because the only difference is the inclusion or exclusion of two additional unchanged values for partition points $1$ and $N$. Indeed, if we chose to include just two unchanged values --- one for each of first and $N$th partition points --- the formulations are equivalent. Thus, the altered $P_B$ may be thought of as a map to \emph{some near point} of the agreement set where the formulation in question is the three set formulation. Because of this similarity, we do not consider these two approaches separately. For all of our examples we use the three set reformulation which does not include unchanged values in the averaging step.

Simple code for computing the projections $P_A$ and $P_B$ may be seen in Algorithm~\ref{alg:PAandPB}. It uses the stored procedures from Algorithm~\ref{alg:generatehypersurfaces}. Note that the projection $P_B$ is the three set reformulation described above which does not include unchanged values in the averaging step.

\begin{algorithm}\caption{Compute Projection for $A = \Omega_1 \times \dots \times \Omega_N$ and a near point in $B$}\label{alg:PAandPB}
	\SetKwFor{proc}{procedure}{}{end:}
	\SetKwFor{uFor}{for}{do}{end:}
	\proc{Project onto A}{
		\KwData{Receives as input a point $x=(x_1,\dots,x_N) \in \mathbb{R}^{N\times N}$ ($x_k \in \mathbb{R}^N$ for all $k$)}
		\KwResult{Returns a point $u=(u_1,\dots,u_N) \in \mathbb{R}^{N\times N}$ such that $u \in P_A(x)$.}
		\uFor{$k \in \{1,\dots,N\}$ }{
			set $u_k:=\textbf{Projection for $\Omega_k$}(u_k)$;
		}
		return $u$;
	}
	
	\proc{Concur in B}{
		\KwData{Receives as input a point $x=(x_1,\dots,x_N) \in \mathbb{R}^{N\times N}$ ($x_k \in \mathbb{R}^N$ for all $k$)}
		\KwResult{Returns a point $u=(\mu,\dots,\mu)\in \mathbb{R}^{N\times N}$ where $\mu \in \mathbb{R}^N$ (Clearly $u \in B$).}
		set  $\mu_1:=\frac{1}{2}(x_{1,1}+x_{2,1})$;\\
		set  $\mu_N:=\frac{1}{2}(x_{{N-1},N}+x_{N,N})$;\\
		\uFor{$j \in \{2,\dots,N-1\}$ }{
			set  $\mu_j:=\frac{1}{3}(x_{j-1,j}+x_{j,j}+x_{j+1,j})$;
		}
		\uFor{$k \in \{1,\dots,N\}$ }{
			set  $u_k:=\mu$;
		}
		return $u$;						
	}
\end{algorithm}

\section{Examples}\label{s:examples}

For all of our examples, unless otherwise specified, we use as an initial point for the iterations $x_0=(\omega,\dots,\omega)\in R^{N\times N}$ where $\omega_i =\alpha+ \frac{i(\beta-\alpha)}{N+1}, i=1,\dots,N$, are the node values of the affine function satisfying the boundary values. We also use $N=21$ unless otherwise specified. We compute the error in the natural way:
\begin{equation}
	\epsilon:= \frac{b-a}{N+1}\sum_{k=1}^{N}|\omega'_k-\omega_k|^2.
\end{equation}
We will use the following terms when discussing the error.
\begin{enumerate}
	\item When $\omega'_k$ is the value of the true solution at $x_k = a+\frac{k(b-a)}{N+1}$ and $\omega_k$ represents the solution of the finite difference problem \eqref{sys} at $x_k$ calculated using Newton's method, $\epsilon$ measures the true error of the approximate solution from the true solution. We expect this error to decrease as $N$ is increased. We show this error for each of our examples with both $N=11$ and $N=21$ in table~\ref{tab:big}.
	
	\item When the $\omega_k$ are values obtained from DR or AP and the $\omega'_k$ are the values at $x_k$ of the true solution, we call $\epsilon$ the \emph{Error from true solution}. 
	
	\item When the $\omega_k$ are values obtained from DR or AP and the $\omega'_k$ are the values at $x_k$ obtained by Newton's method (which is taken to be the numerical solution of the finite difference problem \eqref{sys}) we call $\epsilon$ the \emph{error from Newton solution}.
	
	\item For an iterate of the method of alternating projections (AP) each iterate lies on the agreement set $B$ and so we compute the error where the $\omega_k$ are the induced numerical solution. For each iterate of Douglas-Rachford method (DR) we project back onto $B$ to obtain a numerical solution. In either case, we take \emph{relative error} to mean the change in numerical solution from one iterate to the next: the value of $\epsilon$ when the $w'_k$ and $w_k$ values correspond to numerical solutions from the $n$th and $(n-1)$th steps of the method we are scrutinizing.
\end{enumerate}

\begin{figure}[h]
	\begin{subfigure}{.5\textwidth}
		\begin{center}
			\includegraphics[width=\linewidth]{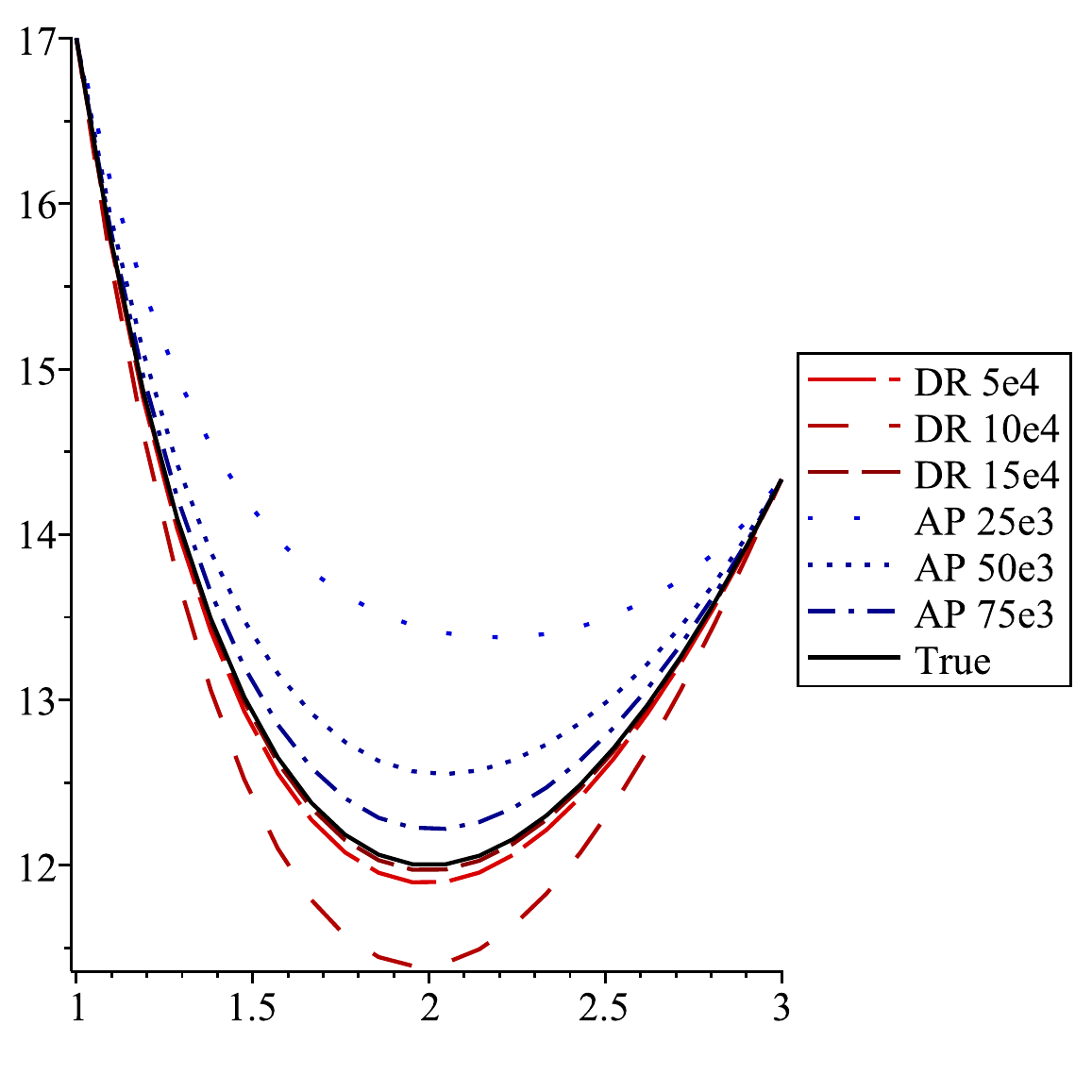}
		\end{center}
		\caption{True and numerical solutions}
	\end{subfigure}
	\begin{subfigure}{.48\textwidth}
		\begin{center}
			\includegraphics[width=\linewidth]{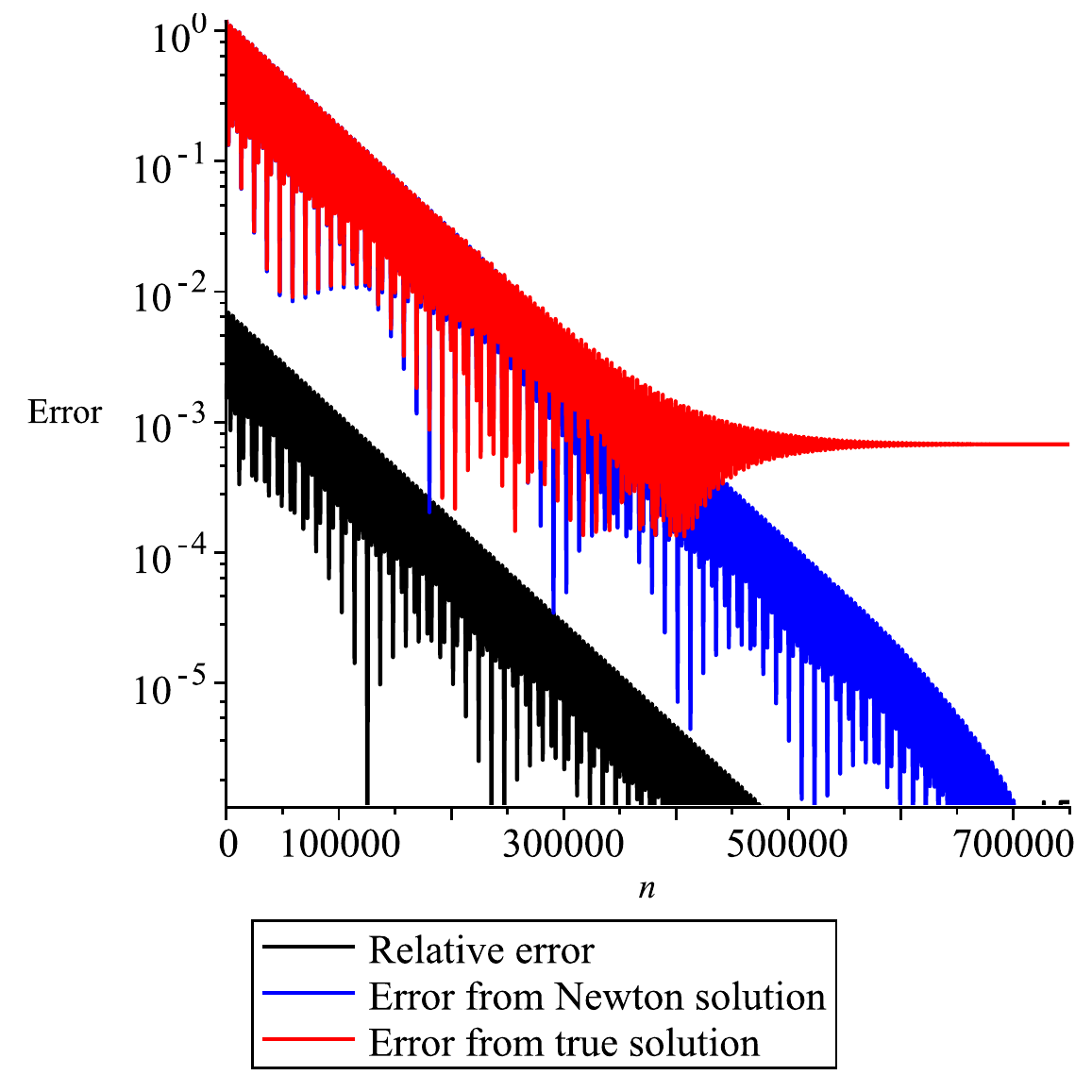}
		\end{center}
		\caption{Error for DR iterates}
	\end{subfigure}
	\caption{Convergence behavior for Example~\ref{ex:Book}} \label{fig:bookcompare}
\end{figure}

When we plot numerical solutions corresponding to various iterates of our methods (as at left in Figure~\ref{fig:bookcompare}), we report first the name of the method (DR or AP) followed by the number of the iterate for which we are plotting a numerical solution. We use the shorthand $N${\scriptsize E}$M := N \cdot 10^M$.

In cases where Newton's method converges, it generally achieves a difference between subsequent iterates of less than $10^{-12}$ within $10$ steps. As will become apparent from the examples, this is so much faster than our methods as to render any comparison of speed useless. However, our methods sometimes work in cases where Newton's method struggles, and they provide useful insights into the behaviour of such algorithms in the nonconvex setting more generally, complementing previous work in this area. The motivating and ideal conditions for implementation are further discussed in the Conclusion (section~\ref{section:conclusion}).

\begin{example}\label{ex:Book}We first tested the method on a simple problem from \cite{BF16}: namely the differential equation $y''=\frac18(32+2x^3-yy')$ with boundary conditions $y(1)=17,y(3)=43/3$, which admits the smooth solution $$y(x)=x^2+16/x.$$ 
\end{example}	
DR, AP, and Newton's method all successfully solve the induced system of equations. Their behaviour is shown graphically in Figure~\ref{fig:bookcompare} where $N=21$.

\begin{figure}[h]
	\begin{subfigure}{.49\textwidth}
		\begin{center}
			\includegraphics[width=\linewidth]{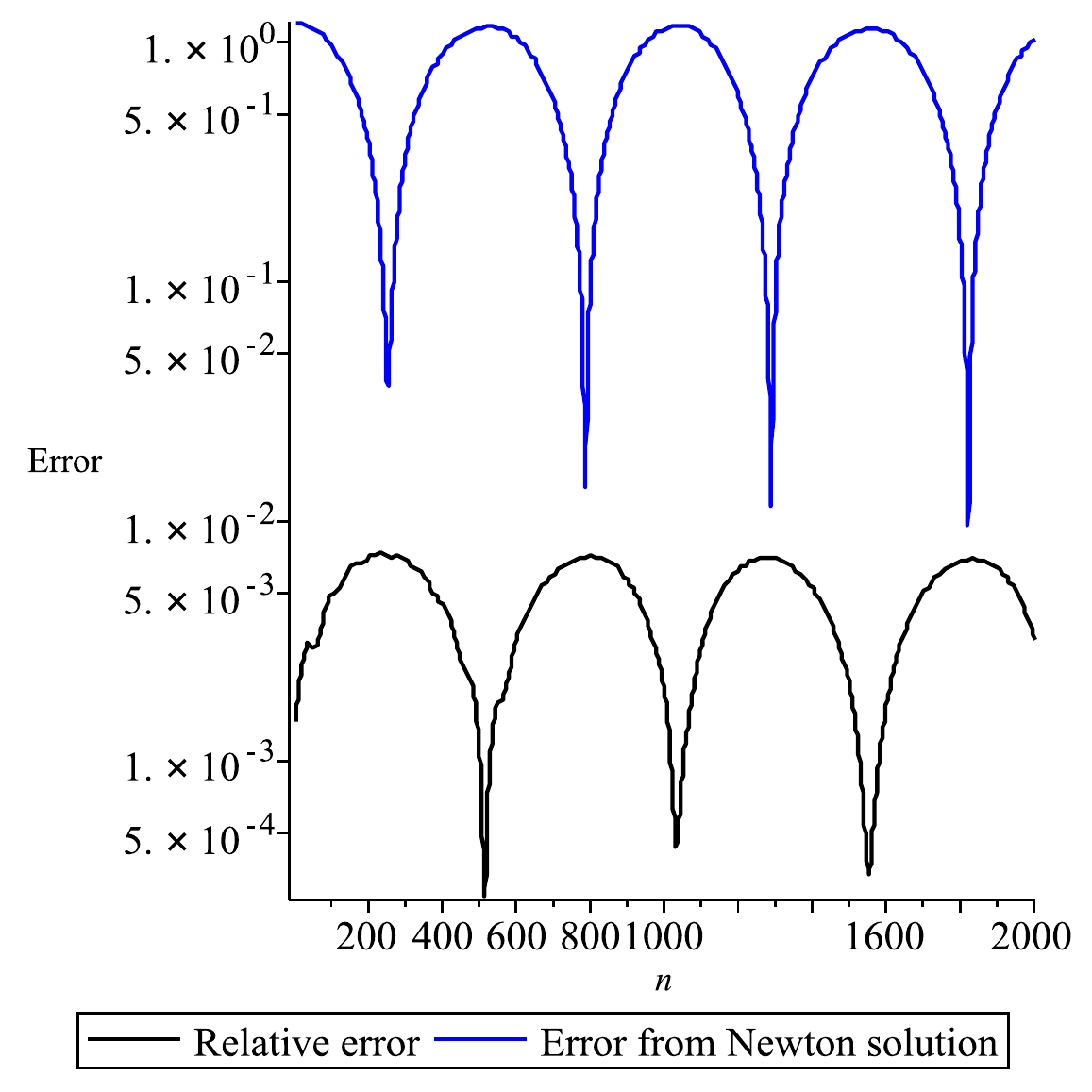}
			\caption{DR}
		\end{center}
	\end{subfigure}
	\begin{subfigure}{.49\textwidth}
		\includegraphics[width=\linewidth]{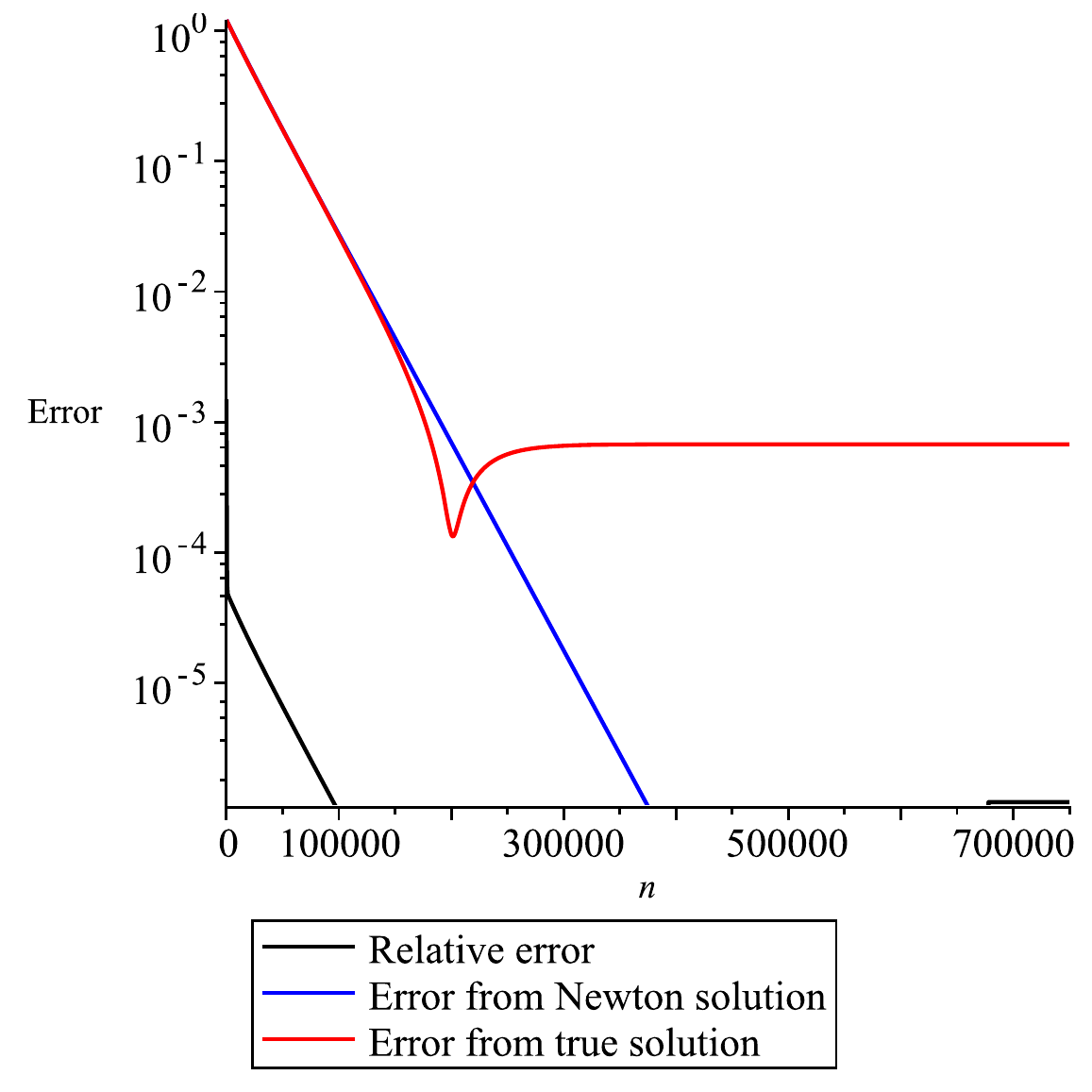}
		\caption{AP}
	\end{subfigure}
	\caption{Convergence behavior for Example~\ref{ex:Book}.} \label{fig:bookerr}
\end{figure}

At around 400,000 iterates, the numerical solution from DR is close to the solution of the finite differences problem \eqref{sys} and so the error from the true solution appears to stabilize, exposing the inherent error between the approximate solution (with $21$ nodes) and the true solution.

Zooming in, the first 2,000 iterates are shown at left in Figure~\ref{fig:bookerr}; we see that the ``solid'' appearance in Figure~\ref{fig:bookcompare} is created by shorter-scale oscillations in relative error. At right in Figure~\ref{fig:bookerr} we see the behavior of AP which converges more quickly and also without the drastic changes in relative error so typical of DR. This pattern of converging faster was observed often though not always, and the relative error plots for AP were similar in all our examples.\bigskip

In the next two examples we consider the effect of partition size on the error from the true solution and on the rate of convergence.

\begin{figure}[h]
	\hspace*{\fill}\includegraphics[width=.49\textwidth]{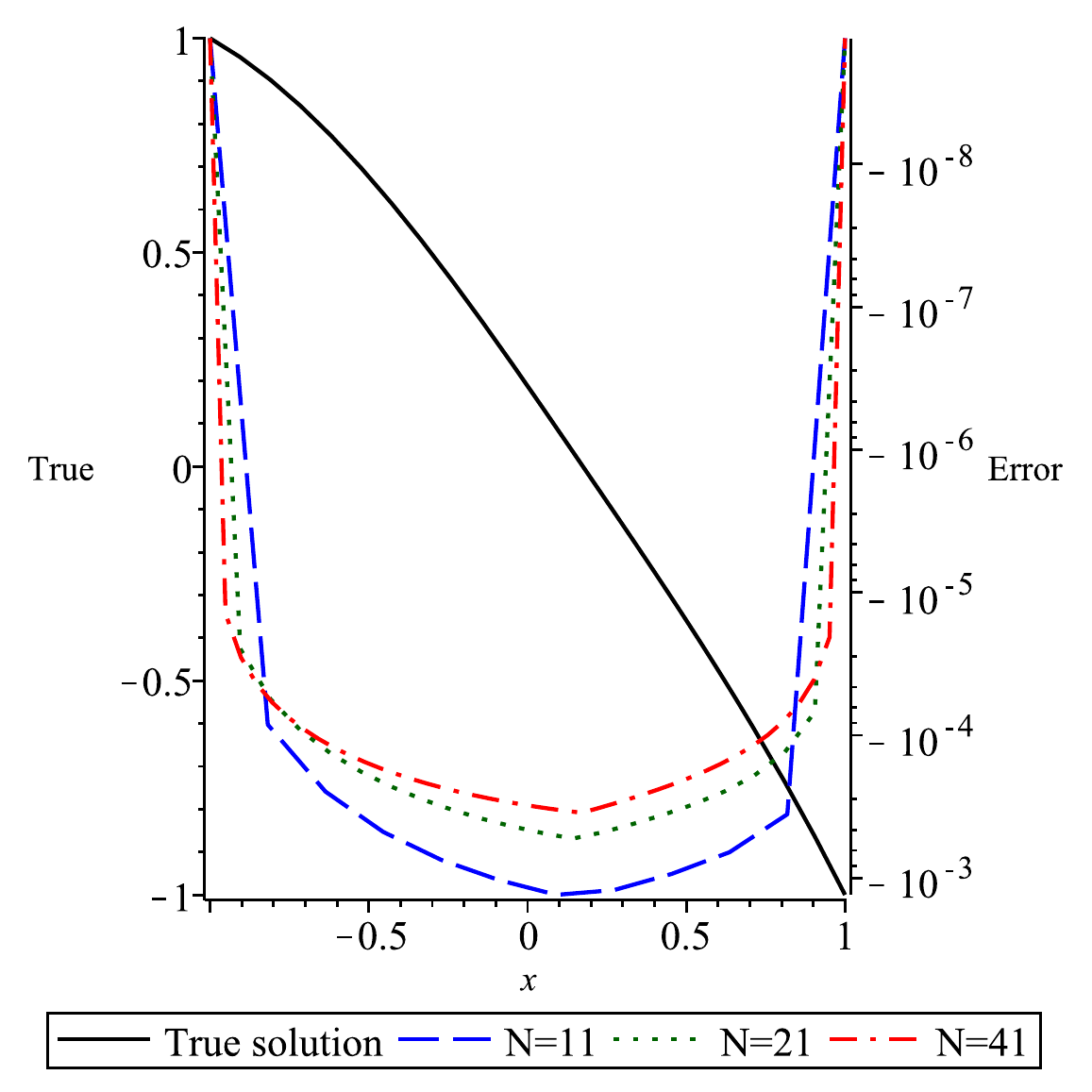}\hspace*{\fill} \includegraphics[width=.49\textwidth]{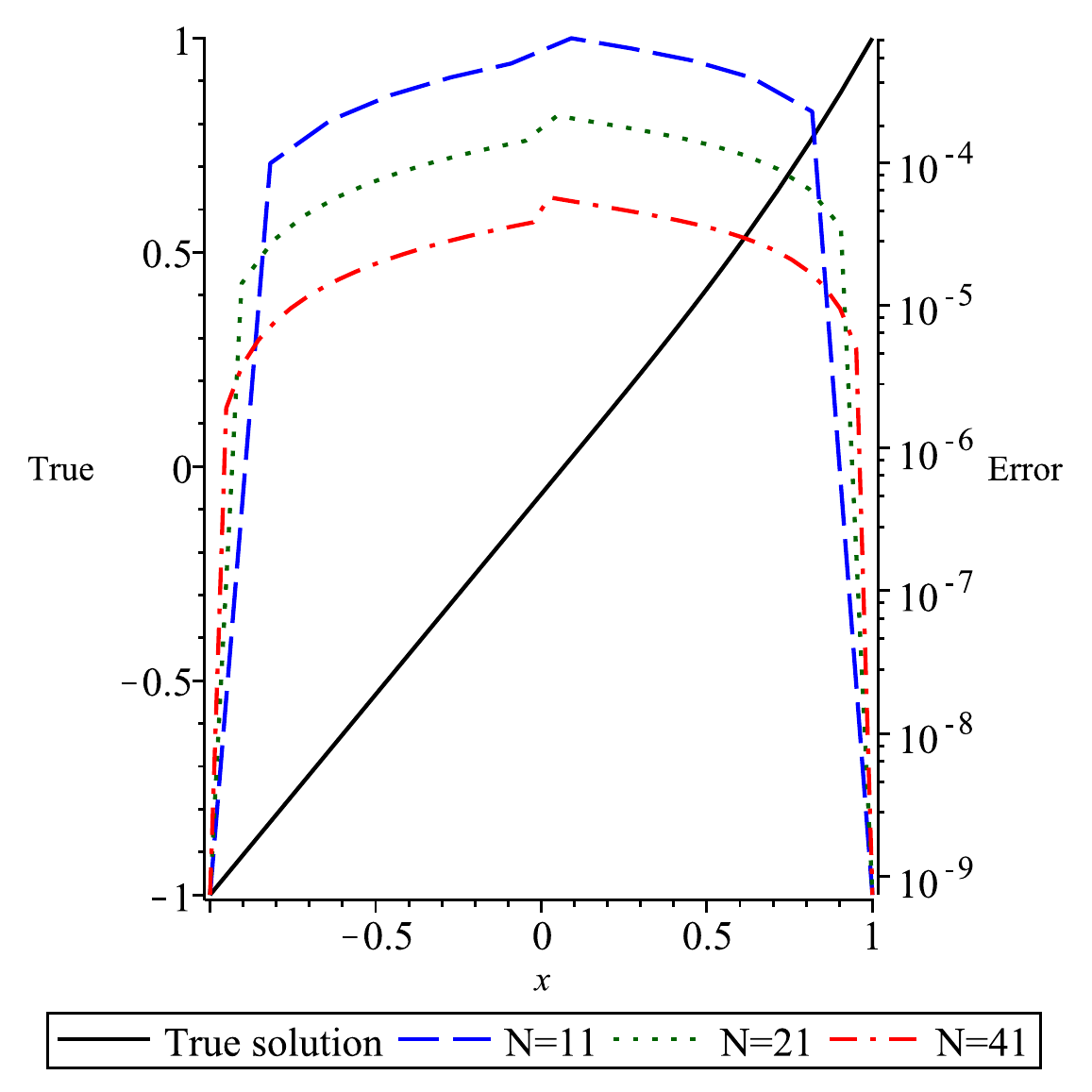}\hspace*{\fill}
	\caption{True solutions (left axis scale) and effect of partition size on error between true solution and estimate by Newton (right axis scale) for Examples~\ref{ex:BraileyAbs} (left) and \ref{ex:Jump} (right).}\label{fig:Neffect}
\end{figure}

\begin{example}\label{ex:BraileyAbs}
	We consider the equation $y''=-|y|$ with boundary conditions $y(-1)=1,y(1)=-1$ which admits the smooth solution
	\begin{align*}
		y(x)&=\begin{cases}
			c_1\sin(x)+c_2\cos(x) & x\leq \frac{1}{2}\log\left(\frac{c_4}{e+c_4}\right)\\
			c_3\exp(x)+c_4\exp(-x) & x>\frac{1}{2}\log\left(\frac{c_4}{e+c_4}\right)
		\end{cases}\\
		c_1&=\frac{c_2 \cos(1)-1}{\sin(1)}\\
		c_2&\resizebox{.9\hsize}{!}{%
			$=\frac{- \left( \tan \left( 1 \right) +\tan \left( \frac{1}{2}\,\log  \left( {\frac{c_{{4}}}{{\rm e}+c_4}} \right)  \right)  \right)}{ \left( \tan \left( 1
				\right) \tan \left( \frac{1}{2} \,\log  \left( {\frac {c_{{4}}}{{\rm e}+c_{{4}}
				}} \right)  \right) \sin \left( 1 \right) -\cos \left( 1 \right) \tan
				\left( 1 \right) -\cos \left( 1 \right) \tan \left( \frac{1}{2}\,\log  \left( 
				{\frac {c_{{4}}}{{\rm e}+c_{{4}}}} \right)  \right) -\sin \left( 1
				\right)  \right)}$%
		}\\
		c_3&=-\frac{c_4 {\rm e}^{-1}+1}{{\rm e}}\\
		c_4&\approx 0.6453425944	.
	\end{align*}
\end{example}	
We found convergence for each of our methods. The true solution and the effect of partition fineness ($N$) on the error between various approximations and the true solution is shown at left in Figure~\ref{fig:Neffect}.\bigskip

\begin{figure}[h]
	\hspace*{\fill}\includegraphics[width=.49\textwidth]{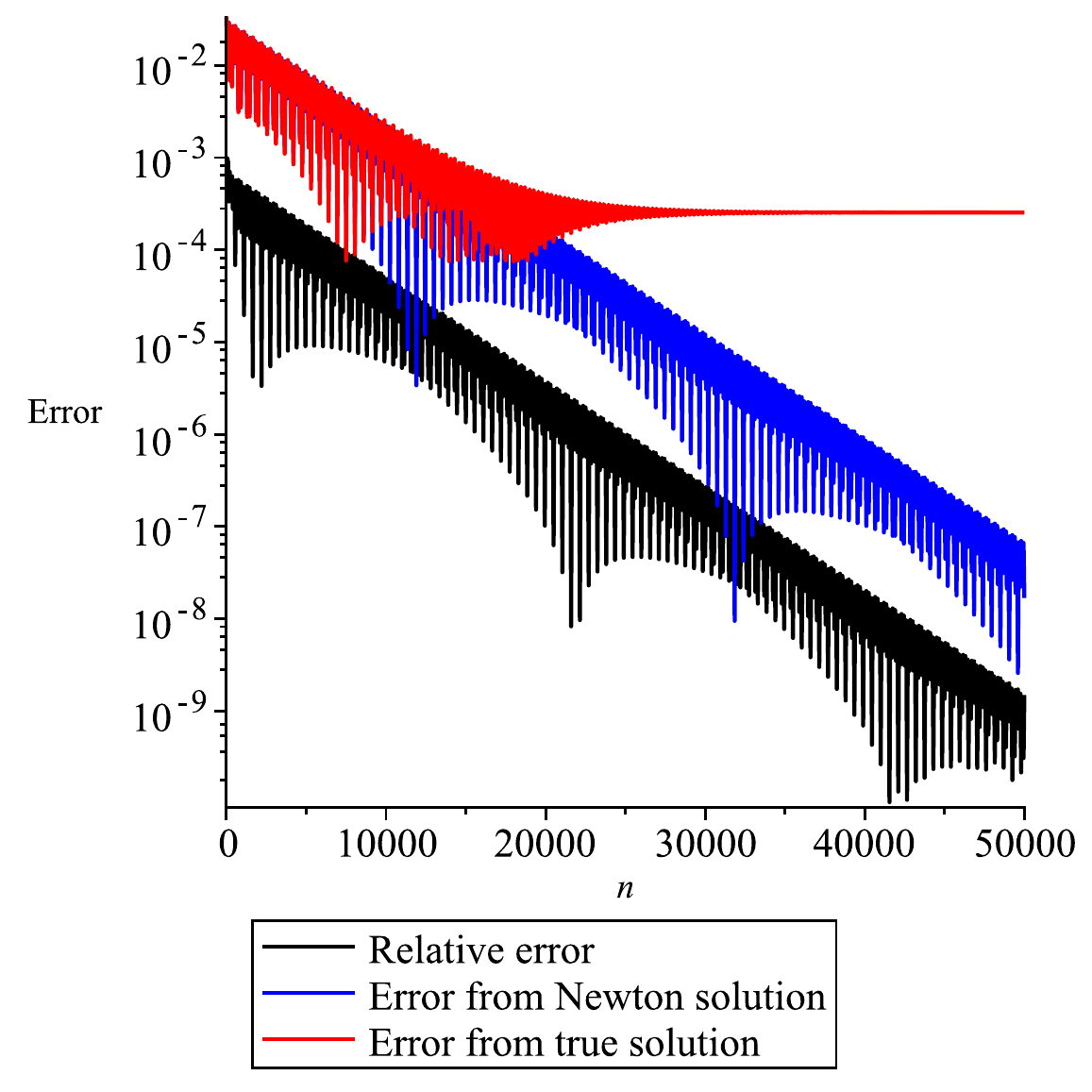}\hspace*{\fill} \includegraphics[width=.49\textwidth]{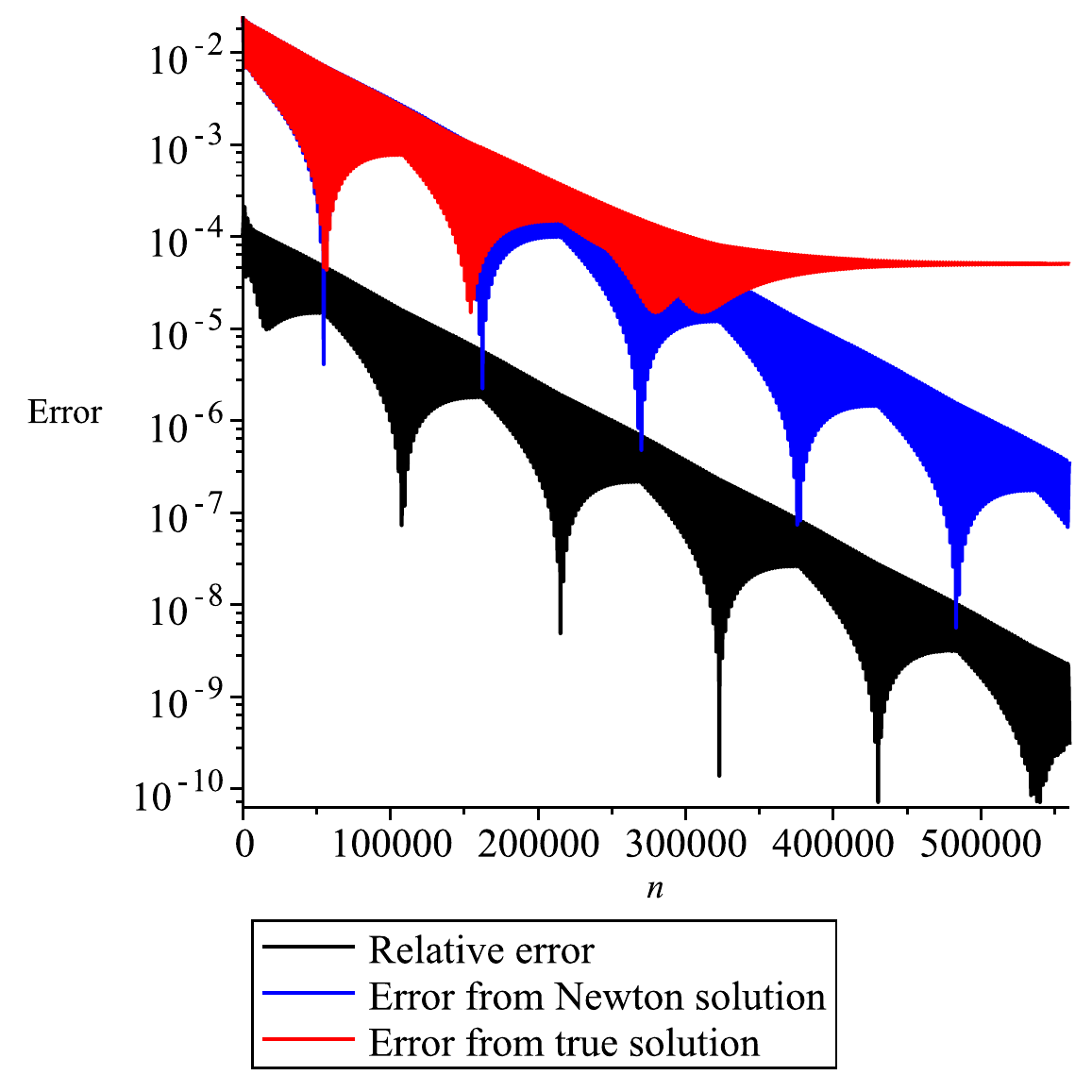}\hspace*{\fill}
	\caption{Effect of $N$ on DR convergence for  Example~\ref{ex:Jump}.}\label{fig:jumpNerr}
\end{figure}

\begin{example}\label{ex:Jump} We examine the differential equation
	\begin{equation}
		y'' = \begin{cases}
			0 & x<0\\
			y & x\geq 0
		\end{cases}
	\end{equation}
	with boundary conditions $y(-1)=-1$ and $y(1)=1$ which admits the smooth solution
	\begin{equation}
		y(x)=\begin{cases}
			\left(\frac{{\rm e}^{-1}+2}{2{\rm e}}+\frac12\right)x+\left(\frac{{\rm e}^{-1}+2}{2{\rm e}}-\frac12 \right) & x < 0\\
			\frac{{\rm e}^{-1}+2}{2{\rm e}}{\rm e}^x-\frac12 {\rm e}^{-x} & x \geq 0.
		\end{cases}
	\end{equation}
	The true solution and the effect of $N$ on the error between true and approximate solutions is shown at right in Figure~\ref{fig:Neffect}. A convergence plot for DR is given in Figure~\ref{fig:jumpNerr} where $N=11$ is shown at left and $N=21$ is shown at right. 
\end{example}	

Noting the different horizontal axis scales, it may be seen that, as one would expect, convergence is much more rapid for smaller $N$, a phenomenon which held both consistently and dramatically across all our examples. 

The ``aqueducts''---which might seem to suggest long-scale oscillations in the change from iterate to iterate---appear to be an artifact of the sample of iterates we used to prepare the plot. For $N=21$ our plot is made from sampling at every $400$th iterate. Shorter scale oscillations of the kind visible in Figure~\ref{fig:bookerr} appeared for all of our error plots, and by sampling infrequently we tend to sample near the tops and sides of the humps while missing the valleys. This phenomenon combined with the regularity of the shorter scale oscillations creates the illusion of aqueducts.

The relative error plots do, however, reveal an important characteristic of the behavior. The change in error from the true solution does \emph{not} track the relative error between iterates but \emph{instead} roughly tracks the change in relative error at the tops of the humps in Figure~\ref{fig:bookerr}. Once sufficiently close to the solution, these oscillations become regular and so convergence can be estimated by tracking only the iterates where relative error peaks.

\begin{figure}
	\begin{center}
		\includegraphics[width=.3\textwidth]{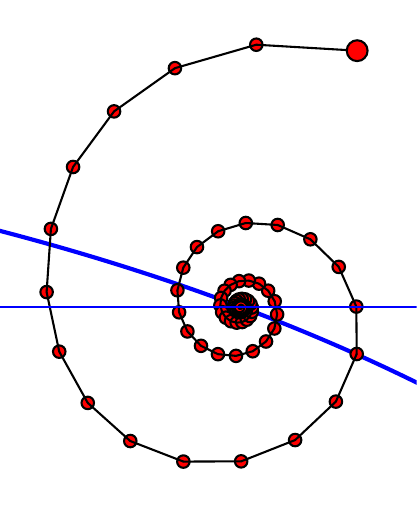} \includegraphics[width=.4\textwidth]{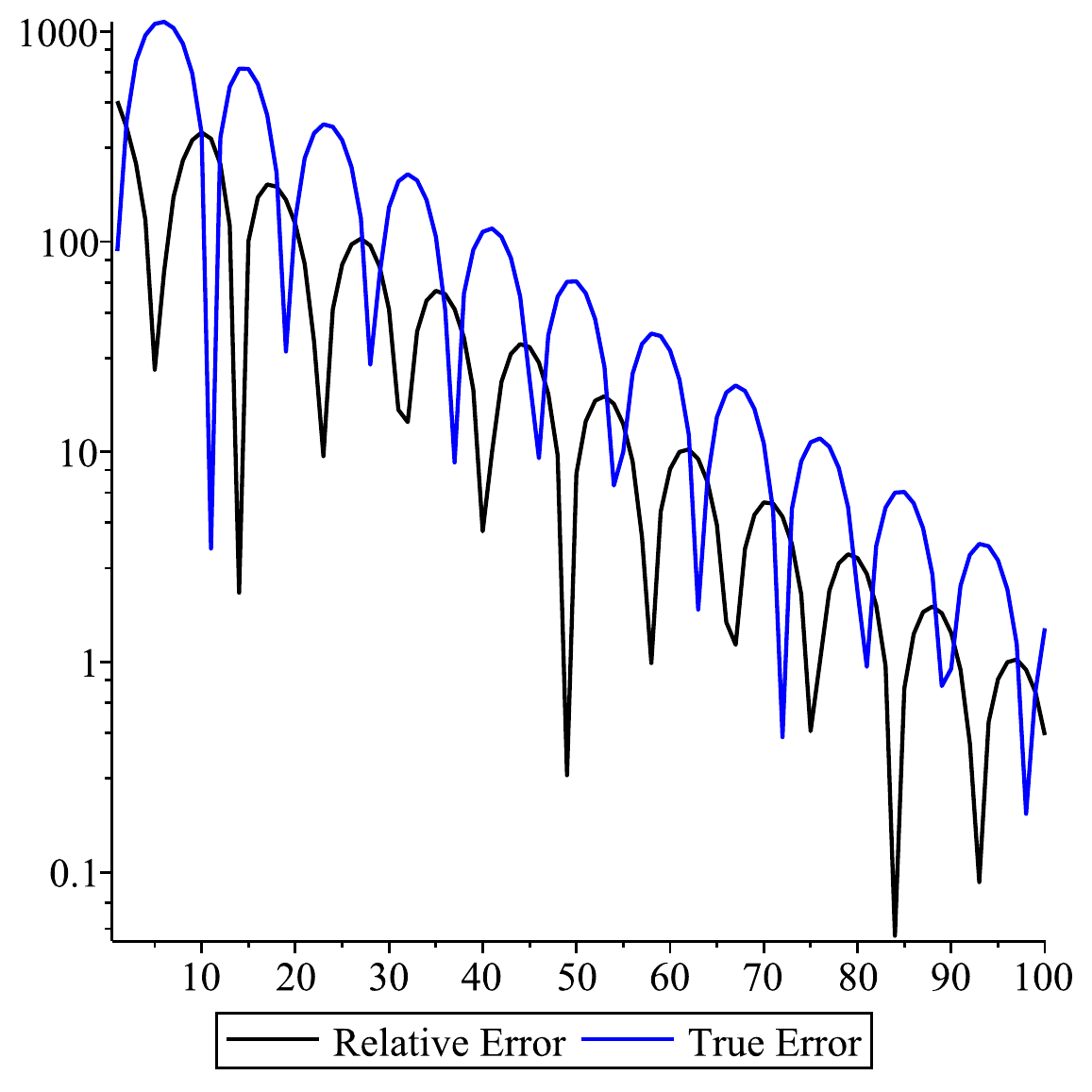}
	\end{center}
	\caption{Relative error and error from true solution for converging DR iterates for an ellipse $E$ and a line $L$.}\label{fig:ellipse2}
\end{figure}

This behavior is consistent with the behavior of DR in other contexts. At left in Figure~\ref{fig:ellipse2} we see DR iterates for an ellipse and a line. The line is the analog of our diagonal set $B$ \eqref{eqn:B}, and so at right we report $\|P_L x_{n+1} - P_L x_{n}\|_2$. The similarities to Figure~\ref{fig:bookerr} are unmistakable.\bigskip

In each of the next three examples we consider the sensitivity of the methods to the starting point. For the first two examples we have multiple potential solutions, and for the final example Newton's method may cycle rather than finding a solution.

\begin{figure}
	\begin{center}
		\includegraphics[width=.49\textwidth]{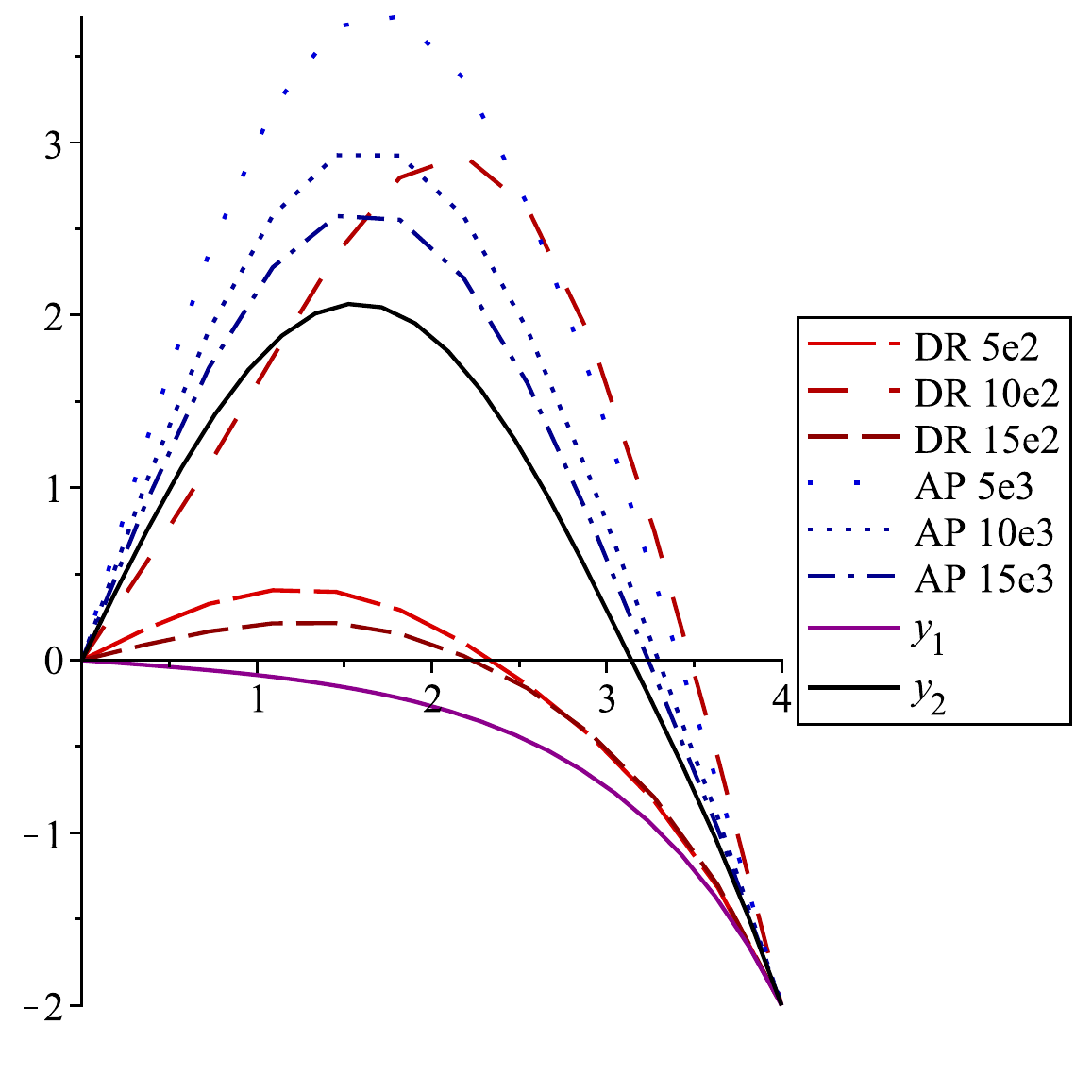} \includegraphics[width=.49\textwidth]{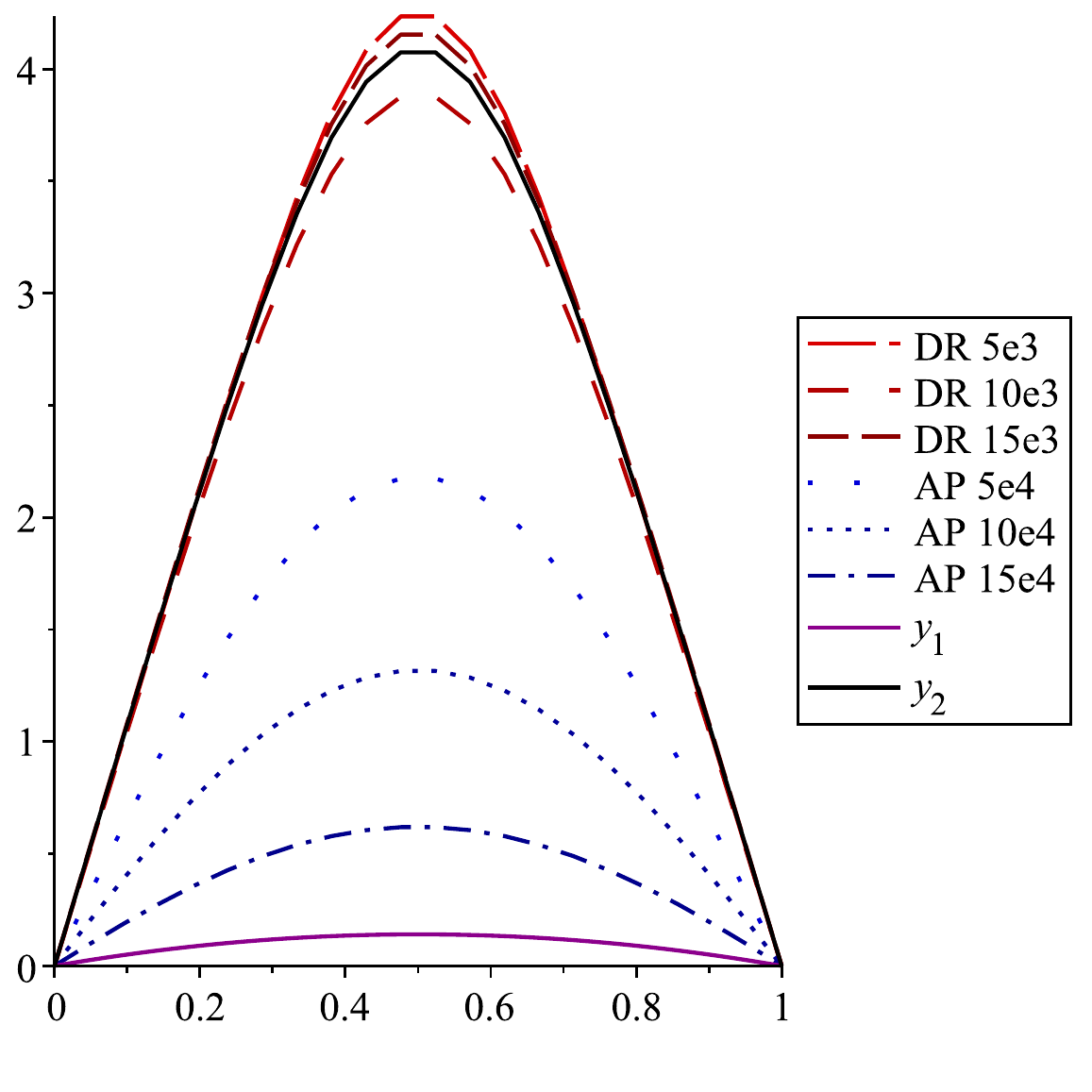}
	\end{center}
	\caption{DR and AP may converge to two different solutions from the same starting point: at left  Example~\ref{ex:PaperAbs}, at right Example~\ref{ex:Exp}.} \label{fig:differentsolutions}
\end{figure}

\begin{example}\label{ex:PaperAbs}The differential equation $ y''=-|y|$ with boundary conditions $y(0)=0, y(4)=-2$ admits two possible smooth solutions:
	\begin{align}
		y_1(x)&=-\frac{2\sinh(x)}{\sinh(4)} \label{eqn:Abs1}\\
		y_2(x)&=\begin{cases}
			\frac {2\sin \left( x \right)}{\sinh(4-\pi)} & x \leq \pi\\
			-\frac{2\sinh(\pi-x)}{\sinh(\pi-4)} & x>\pi \end{cases}.\label{eqn:Abs2}
	\end{align}	
\end{example}	
Here we intially found convergence for small $N$, but our scripts stopped working for larger $N$. Investigating, we found that \emph{Maple}'s \emph{fsolve} was unable to compute the solution to the Lagrangian system for the $P_{\Omega_i}$. Replacing it with our own numerical solver, we recovered convergence. With the starting values corresponding to the affine function matching the boundary conditions, all methods converged to solution $y_1$ from \eqref{eqn:Abs1}. However, with the starting values matching the boundary conditions and $4\chi_{(0,1)}$ everywhere else, AP goes to the ``nearer'' solution of $y_2$ \eqref{eqn:Abs2} while DR finds its way down to $y_1$. This may be seen in Figure~\ref{fig:differentsolutions}. 	

We repeated the experiment for a variety of starting points corresponding to functions which matched the boundary values and were $\lambda\chi_{(0,4)}$ everywhere else for various $\lambda$. The results are tabulated in Table~\ref{tab:abs}. 

\begin{table}
	\begin{tabular}{l | c c c c c c c c c c c c}
		Method/Start $\lambda$ & .01 & .1 & .5 & 1 & 2 & 3 & 4 & 5 & 6 & 7 & 8 & 9\\ \hline
		Newton N=11 & $2$ & $2$ & $2$ & $2$ & $2$ & $2$ & $2$ & $2$ & $2$ & $2$ & $2$ & $2$\\
		DR N=11 & $1$ & $1$ & $2$ & $2$& $2$ & $2$ & $1$ & $1$ & $1$ & $1$ & $1$ & $1$ \\
		AP N=11 & $1$ & $1$ & $2$& $2$& $2$& $2$& $2$& $2$& $2$& $2$ & $2$ & $2$ \\ \hline
		Newton N=21 & $2$ & $2$ & $2$ & $2$ & $2$ & $2$ & $2$ & $2$ & $2$ & $2$ & $2$ & $2$ \\
		DR N=21 & $1$ & $1$ & $2$ & $2$ & $2$ & S & S & S & S & $2$ & $2$ & $2$ \\
		AP N=21 & $1$ & $1$ & $2$ & $2$ & $2$ & $2$ & $2$ & $2$ & $2$ & $2$ & $2$ & $2$ \\ \hline
		
		Method/Start $\lambda$ & -.01 & -.1 & -.5 & -1 & -2 & -3 & -4 & -5 & -6 & -7 & -8 & -9\\ \hline
		Newton N=11 & $1$ & $1$ & $1$ & $1$ & $1$ & $1$ & $1$ & $1$ & $1$ & $1$ & $1$ & $1$ \\
		DR N=11 & $1$ & $1$ & $1$ & $1$ & $1$ & $1$ & $1$ & $1$ & $1$ & $1$ & $1$ & $1$\\
		AP N=11 & $1$ & $1$ & $1$ & $1$ & $1$ & $1$ & $1$ & $1$ & $1$ & $1$ & $1$ & $1$ \\ \hline
		Newton N=21 & $1$ & $1$ & $1$ & $1$ & $1$ & $1$ & $1$ & $1$ & $1$ & $1$ & $1$ & $1$ \\
		DR N=21 & $1$ & $1$ & $1$ & $1$ & $1$ & S & S & S & S & $2$ & $2$ & $2$ \\
		AP N=21 & $1$ & $1$ & $1$ & $1$ & $1$ & $1$ & $1$ & $1$ & $1$ & $1$ & $1$ & $1$
	\end{tabular}
	\caption{Sensitivity to starting point for Example~\ref{ex:PaperAbs}: $1$ or $2$ indicate the method converged to $y_1$ or $y_2$ while S indicates the method appeared stuck after $5${\scriptsize E}$5$ iterates.}\label{tab:abs}
\end{table}

Newton's method behaved very predictably, always converging to $y_1$ for $\lambda <0$ and $y_2$ for $\lambda >0$ regardless of partition size. AP was slightly less predictable, converging to $y_1$ for $\lambda=0.01$. For $\lambda=0.1$ it appeared stuck between $y_1$ and $y_2$ even after $15${\scriptsize E}$4$ iterates regardless of partition size; eventually it converged to $y_1$. 

\begin{figure}
	\begin{center}
		\vspace*{\fill}\includegraphics[width=.5\textwidth]{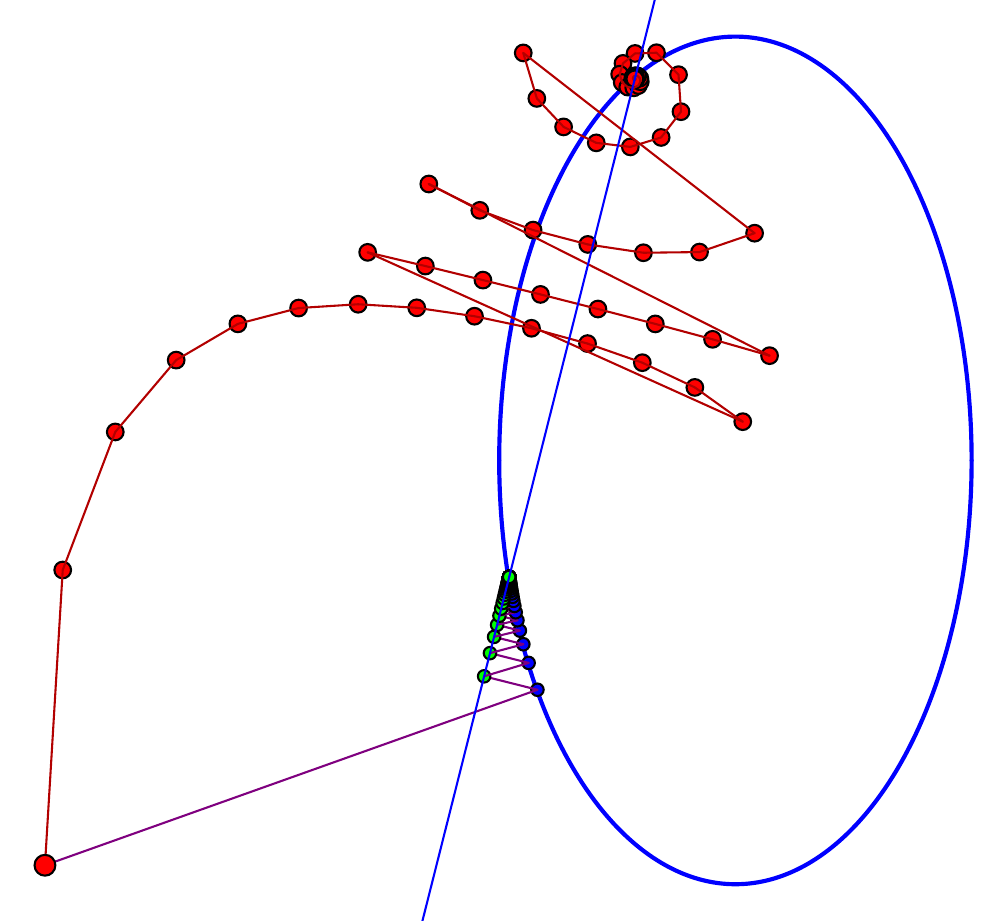}\vspace*{\fill}\includegraphics[width=.48\textwidth]{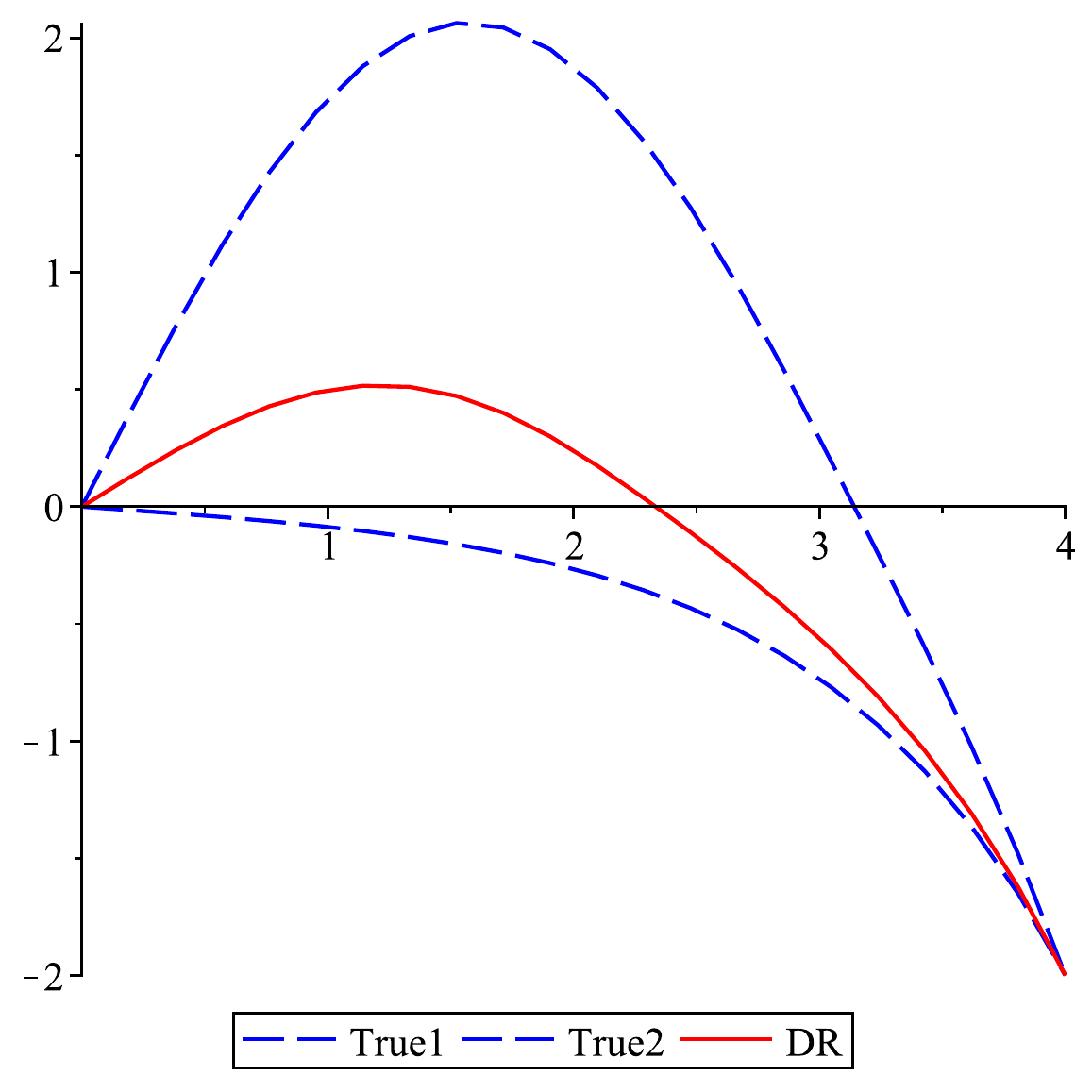}\vspace*{\fill}
	\end{center}
	\caption{Left: DR started sufficiently far from two feasible points may converge to the farther of the two while AP converges to the nearer. Right: for Example~\ref{ex:PaperAbs} after $5${\scriptsize E}$5$ iterates DR appears stuck for some starting points.}\label{fig:ellipse1}
\end{figure}

The behavior of DR, by contrast, was highly unpredictable, changed drastically with partition size, and frequently converged to the ``farther'' away of the two solutions when started some distance from both. This is consistent with the known behavior of Douglas-Rachford illustrated in Figure~\ref{fig:ellipse1}. See, for example, \cite{BLSSS}.

We observed another trend as well. For $\lambda=4$ and $N=11$ DR converged to $y_1$ while for $\lambda=2$ it converged to $y_2$; for $\lambda=3$ convergence was extremely delayed. For most values, we were able to ascertain the eventual solution within $15${\scriptsize E}$4$ iterates. For some $\lambda$ values we were unable to tell even after $5${\scriptsize E}$5$ iterates. This pattern of ``crossroad'' points taking longer to close on a destination held consistently. One example is shown at right in Figure~\ref{fig:ellipse1}. \bigskip

\begin{example}\label{ex:Exp}The differential equation $y'' =-\exp(y)$ with boundary conditions $y(0)=y(1)=0$ admits two smooth solutions:
	\begin{align}
		y(x)&=\log\left(c-c\tanh^2\left(\sqrt{\frac{c}{2}}(1/2-x) \right) \right) \nonumber \\
		\text{where}\; &c\approx 1.1508 \label{expsol1} \\
		\text{or}\; &c \approx 59.827 \label{expsol2}
	\end{align}	
\end{example}	
When the starting values match the unique affine function corresponding to the boundary conditions, all of the numerical methods converged to the particular solution given by \eqref{expsol1} which we call $y_1$. If we start instead from a function matching the boundary conditions and $2\chi_{(0,1)}$ everywhere else, for $N=21$ AP still goes to $y_1$ while DR converges to the other solution $y_2$ given by \eqref{expsol2}. This can be seen in Figure~\ref{fig:differentsolutions}.

We again repeated the experiment for a variety of starting points corresponding to functions which matched the boundary values and were $\lambda\chi_{(0,1)}$ elsewhere for various $\lambda$. The results are tabulated in Table~\ref{tab:exp} where it may be seen that for certain starting values Newton's method diverged or AP appeared stuck after $15${\scriptsize E}$4$ iterates. \bigskip

\begin{table}
	\begin{center}
		\begin{tabular}{l | c c c c c c c c c}
			Method / Start $\lambda$ & -1 & 0 & 1 & 2 & 3 & 4 & 5 & 6 & 7\\ \hline
			Newton N=11 & $1$ & $1$ &$1$ &$1$ &$2$ & D & D & D & D\\
			DR N=11 & $1$ & $1$ & $1$ & $2$ & $2$ & $2$ & $2$ & $2$ & $2$ \\
			AP N=11 & $1$ & $1$ & $1$ & $1$ & $2$ & $2$ & $2$ & S & S \\ \hline
			Newton N=21 & $1$ & $1$ & $1$ & $1$ & $2$ & D & D & D & D\\
			DR N=21 & $1$ & $1$ & $1$ & $2$ & $2$ & $2$ & $2$ & $2$ & $2$\\
			AP N=21 & $1$ & $1$ & $1$ & $1$ & $2$ &  $2$ &$2$ & $2$ & S
		\end{tabular}
	\end{center}
	\caption{Sensitivity to starting point for Example~\ref{ex:Exp}: $1,2$ indicates the method converged to $y_1,y_2$ respectively while ``D'' and ``S'' respectively indicate the method diverged or appeared to hover.}\label{tab:exp}
\end{table}

\begin{figure}
	\begin{center}
		\includegraphics[width=.49\textwidth]{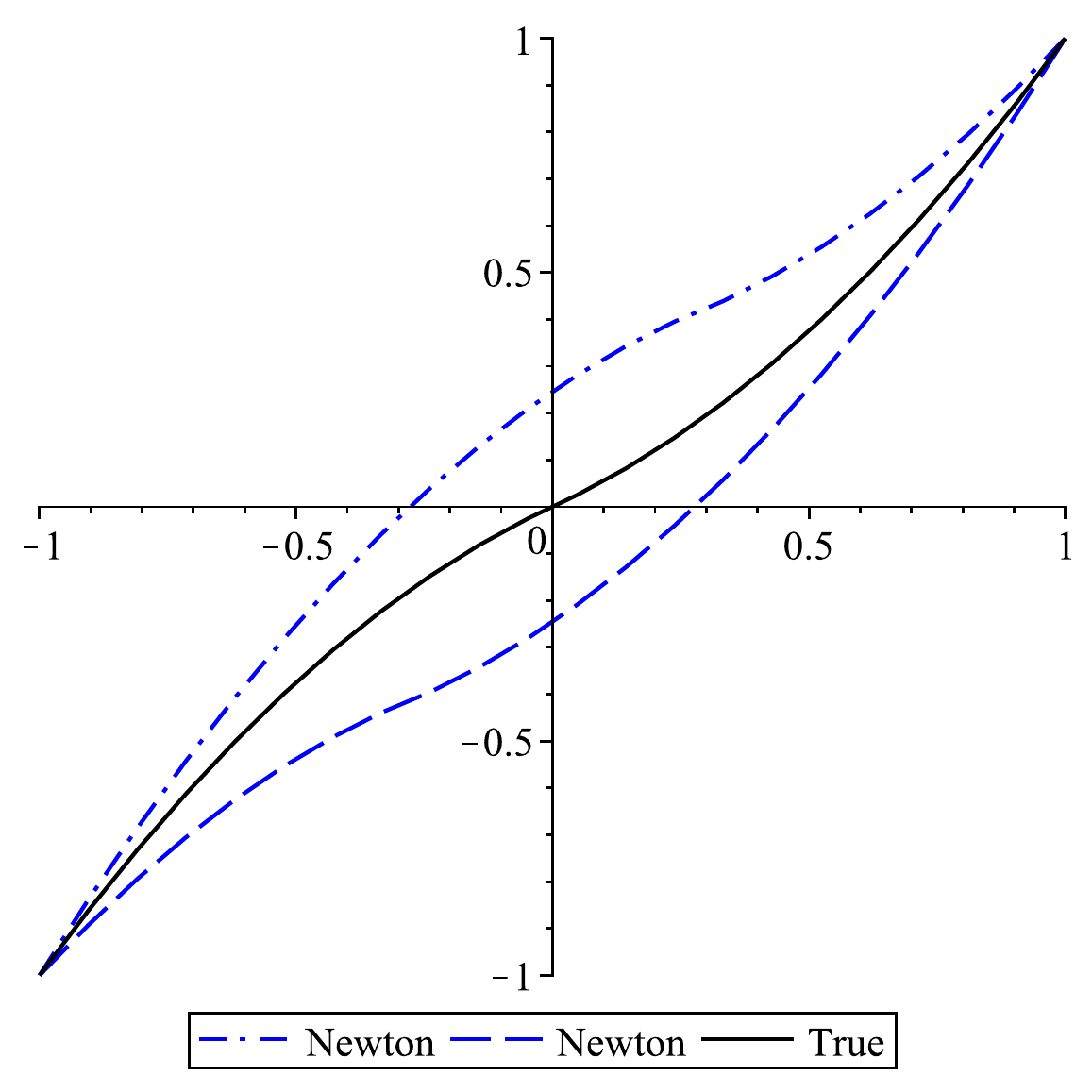}\includegraphics[width=.49\textwidth]{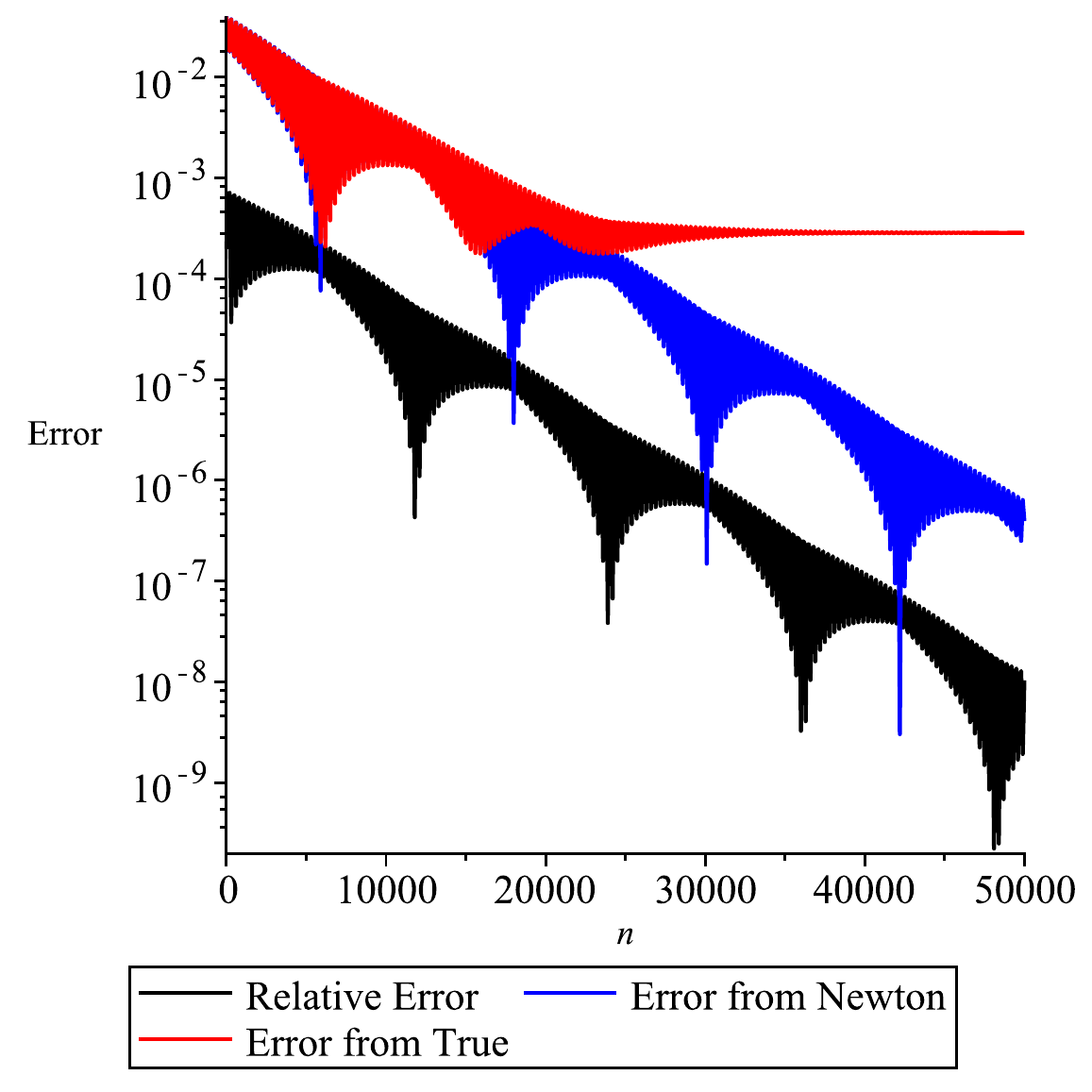}
	\end{center}
	\caption{Newton's Method may cycle for certain starting points in Example~\ref{ex:Heaviside} (left) while DR converges (right).} \label{fig:Heaviside}
\end{figure}

\begin{example}\label{ex:Heaviside} We consider the second order differential equation
	\begin{equation}
		y''(x) = \begin{cases}
			-1 & y(x)<0\\
			1 & y(x)\geq 0,
		\end{cases}
	\end{equation}
	together with the boundary conditions $y(-1)=-1$ and $y(1)=1$.
	Here the right hand side, being a Heaviside function, fails to satisfy the standard conditions for existence and uniqueness. Nonetheless it is readily seen to admit a unique continuous solution on the interval $[-1,1]$, namely the odd function:
	\begin{equation}
		y(x) = \begin{cases}
			-\frac{1}{2}x^2+\frac12 x & x < 0\\
			\frac12 x^2+\frac12 x & x\geq 0.
		\end{cases}
	\end{equation}
\end{example}
This example is especially interesting because while Newton's method finds the solution when starting from the affine function satisfying the boundary criteria, it fails to converge to the solution when started at $1\chi_{(-1,1)},\dots,7\chi_{(-1,1)}$. Instead it cycles between the two non solutions shown at left in Figure~\ref{fig:Heaviside}.

Within $6$ iterates of Newton's method, the norm of the difference between subsequent even iterates or subsequent odd iterates is less than $1${\scriptsize E}$-19$. By way of contrast, DR and AP appear to work well from all of these starting points. At right we show a plot of relative error for DR with $21$ iterates starting from the affine function values.

\begin{table}
	\begin{tabular}{r | r r r r r r}
		& DR & AP & DR & DR & AP & True \\ 
		& $1${\scriptsize E}$-1$ & $1${\scriptsize E}$-1$ & wave & $\frac{{\rm Error}}{{\rm Relative}}$ & $\frac{{\rm Error}}{{\rm Relative}}$ & error\\ \hline
		Ex~\ref{ex:Book}	N=11 &	$9${\scriptsize E}$3$ & $4${\scriptsize E}$3$ & $142$ & $44$ & $2${\scriptsize E}$3$ & $3.4${\scriptsize E}$-3$\\
		N=21 &	$129${\scriptsize E}$3$ & $60${\scriptsize E}$3$ & $516$ & $155$ & $26${\scriptsize E}$3$ & $6.7${\scriptsize E}$-4$\\ \hline
		Ex~\ref{ex:BraileyAbs}	N=11 &	$18${\scriptsize E}$3$ & $9${\scriptsize E}$3$ & $198$ & $63$ & $4${\scriptsize E}$3$ & $4.7${\scriptsize E}$-4$\\
		N=21 &	$247${\scriptsize E}$3$ & $102${\scriptsize E}$3$ & $715$ & $227$ & $53${\scriptsize E}$3$ & $1.3${\scriptsize E}$-4$ \\ \hline
		Ex~\ref{ex:Jump}	N=11 &	$9${\scriptsize E}$3$ & $4${\scriptsize E}$3$ & $138$ & $43$ & $2${\scriptsize E}$3$ & $2.5${\scriptsize E}$-4$\\
		N=21 &	$117${\scriptsize E}$3$ & $58${\scriptsize E}$3$ & $500$ & $155$ & $25${\scriptsize E}$3$ & $5.1${\scriptsize E}$-5$\\ \hline
		Ex~\ref{ex:PaperAbs}	N=11 &	$2${\scriptsize E}$3$ & $1${\scriptsize E}$3$ & $65$ & $19$ & $4${\scriptsize E}$2$ & $3.1${\scriptsize E}$-3$\\
		N=21 &	$25${\scriptsize E}$3$ & $12${\scriptsize E}$3$ & $230$ & $67$ & $5${\scriptsize E}$3$ & $6.2${\scriptsize E}$-4$\\ \hline
		Ex~\ref{ex:Exp}	N=11 &	$16${\scriptsize E}$3$ & $8${\scriptsize E}$3$ & $184$ & $57$ & $34${\scriptsize E}$2$ & $2.6${\scriptsize E}$-5$\\
		N=21 &	$208${\scriptsize E}$3$ & $104${\scriptsize E}$3$ & $670$ & $211$ & $46${\scriptsize E}$3$ & $5.1${\scriptsize E}$-6$\\ \hline
		Ex~\ref{ex:Heaviside}	N=11 &	$1${\scriptsize E}$3$ & $4${\scriptsize E}$2$ & $41$ & $12$ & $1${\scriptsize E}$2$ & $1.4${\scriptsize E}$-3$ \\
		N=21 &	$11${\scriptsize E}$3$ & $5${\scriptsize E}$3$ & $149$ & $46$ & $2${\scriptsize E}$3$ & $2.9${\scriptsize E}$-4$\\ \hline	
	\end{tabular}
	\caption{A summary of experimental results from all examples.}\label{tab:big}
\end{table}

We provide an overview summary of our experimental results in Table~\ref{tab:big}. In the first column we report how many iterates it took for $\log_{10}$ of the ``peak'' relative error for DR to go down by $1$. In the second column we report this for AP where peaks need no longer be considered. In the third column we give the average number of iterates which compose the individual oscillations in the relative error of DR (as in Figure~\ref{fig:bookerr}). In the fourth column we report for DR the ratio of peak \emph{error from the approximate solution} to peak \emph{relative error}. Because the two different peaks do not coincide, we take each peak in the error from the approximate solution and compare it to the \emph{previous} peak in the relative error. In the fifth column we report for AP the ratio of error from approximate solution to the relative error; in this case peaks no longer need be considered. In the final column we show the error of the approximate solution \eqref{sys} (obtained by Newton's method) from the true solution.


\begin{figure}[h]
	\hspace*{\fill}\includegraphics[width=.49\textwidth]{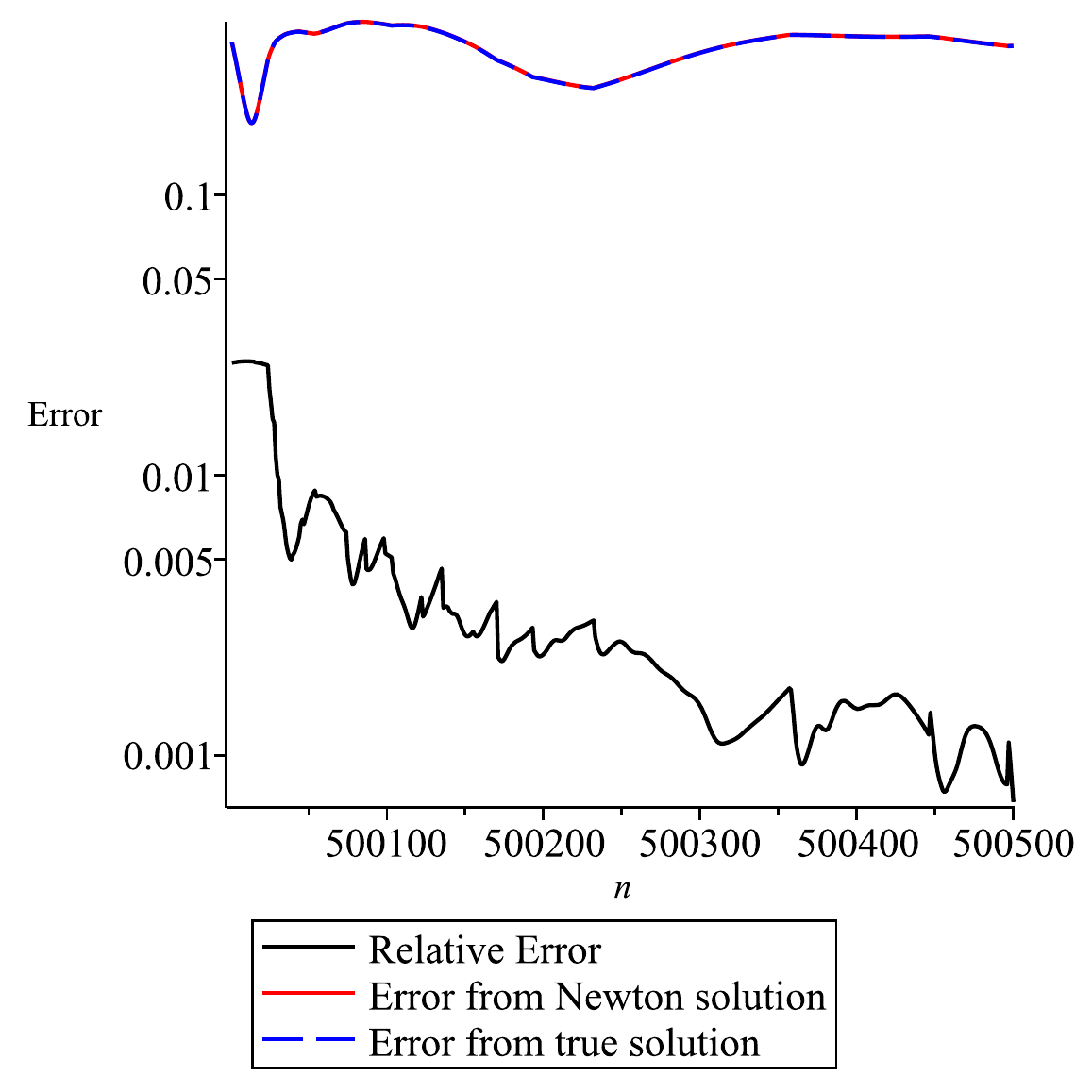}\hspace*{\fill} \includegraphics[width=.49\textwidth]{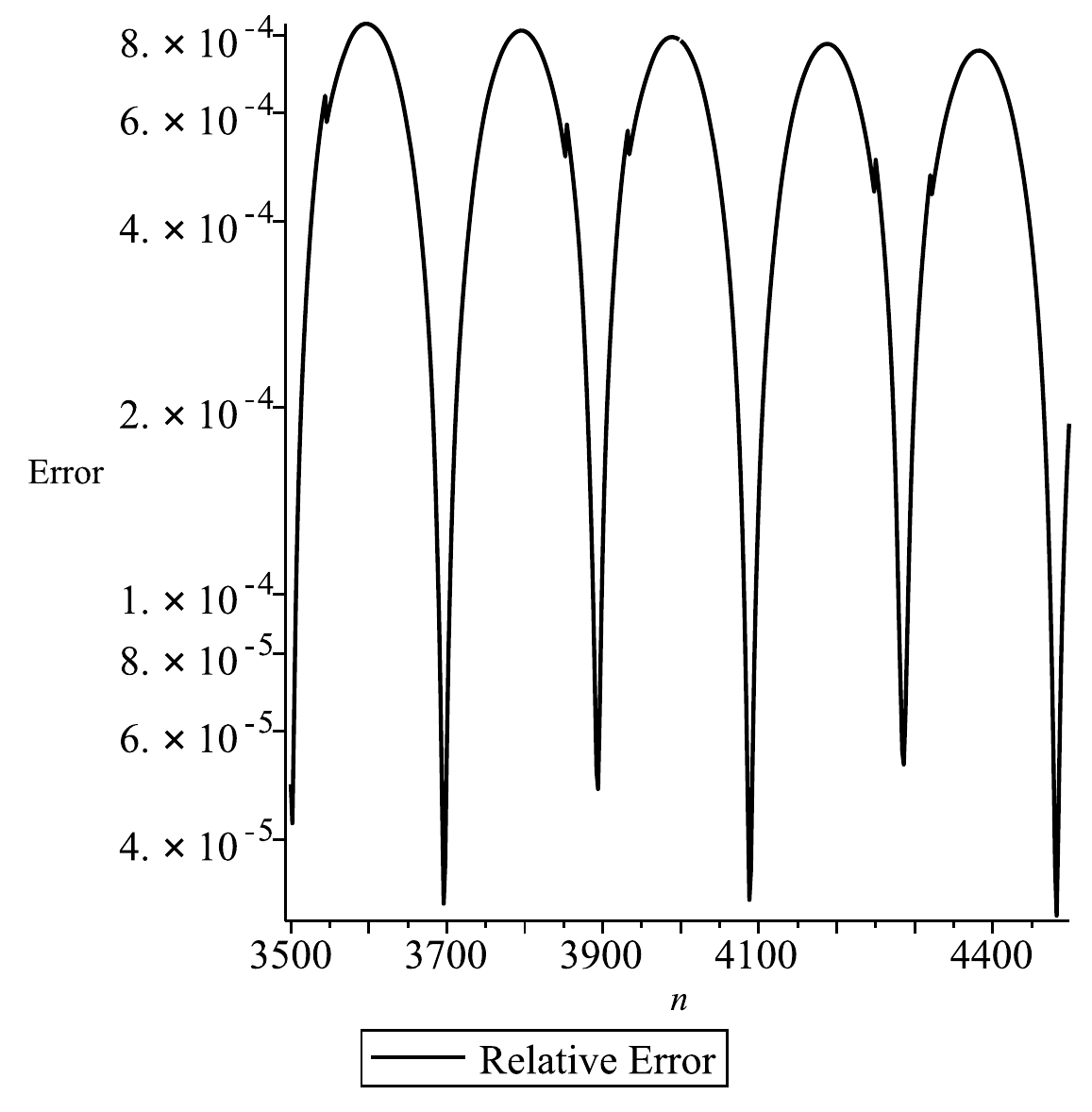}\hspace*{\fill}
	\caption{Left: stuck DR. Right: relative error tends toward a pattern other than smooth oscillation.}\label{fig:relativeerror}
\end{figure}

Analysis of a stuck problem revealed that regular oscillations in relative error were conspicuously absent. This is shown at left in Figure~\ref{fig:relativeerror} where for Example~\ref{ex:PaperAbs} with $N=21$ and starting with $\lambda=6$, DR appears stuck after $5${\scriptsize E}$5$ iterates. Original attempts to catalogue average oscillation length for relative error resulted in data which appeared at times periodic. This led to the discovery that the pattern in relative error may tend toward a predictable pattern other than smooth oscillation. This is shown for Example~\ref{ex:BraileyAbs} with $N=11$ at right in Figure~\ref{fig:relativeerror}.


\section{Conclusion}\label{section:conclusion}

The poor tradeoff in convergence rate for finer partitions suggests some modifications to the method for solving real world problems. One such modification is to begin with a coarse partition and increase the fineness over time. Another is to simply switch to a more traditional solver once sufficient proximity to the true solution is suspected from analysis of the relative error.

The impressive stability of the Douglas-Rachford method relative to more traditional methods is consistent with previous findings in the application of these methods to finding the intersections of analytic curves \cite{LSS}. This property and its unique suitability for parallelization make it an ideal candidate for employment in settings where traditional solvers fail, or for getting close enough to a solution that traditional solvers may be applied.

\vspace{0.5cm}
\textbf{Acknowledgement.} The authors wish to thank an anonymous referee for their careful reading and detailed, helpful feedback.

\end{document}